\documentclass[12pt,reqno]{amsart}
\usepackage{latexsym,amsmath,amsthm,amssymb,amsfonts}
\usepackage{setspace}
\usepackage{color}
\usepackage{tikz}
\usetikzlibrary{patterns}

\numberwithin{equation}{section}
\newtheorem{theorem}{Theorem}[section]
\newtheorem{lem}[theorem]{Lemma}
\newtheorem{thm}[theorem]{Theorem}
\newtheorem{pro}[theorem]{Proposition}
\newtheorem{cor}[theorem]{Corollary}

\newtheorem{rem}[theorem]{Remark}
\newtheorem{eg}[theorem]{Example}

\def\s{\,\,\,\,}

\def\dint{\displaystyle{\int}}

\def\endproof{$\hfill\Box$\\}
\def\R{\mathbb{R}}

\def\vol{\mathrm{vol}}
\def\diam{\mathrm{diam}}
\allowdisplaybreaks

\title[Manifolds for which Huber's Theorem holds]{\bf Manifolds for which Huber's Theorem holds}

\author{Yuxiang Li, Zihao Wang}

\address{\newline
Yuxiang Li:
Department of Mathematical Sciences, Tsinghua University, Beijing 100084, P.R. China.
{\tt Email:liyuxiang@tsinghua.edu.cn}
\newline
\newline
Zihao Wang:
Department of Mathematical Sciences, Tsinghua University, Beijing 100084, P.R. China.
{\tt Email:wangziha19@mails.tsinghua.edu.cn}}

\thanks{The first author is partially supported by NSFC, Grant No.12141103
}
\date{}
\begin{document}
\maketitle

\begin{abstract}
Extensions of Huber's Theorem to higher dimensions with $L^\frac{n}{2}$ bounded scalar curvature
have been extensively studied over the years. In this paper, we 
delve into the properties of  conformal metrics  on a punctured ball with
$\|R\|_{L^\frac{n}{2}}<+\infty$,  aiming to identify necessary geometric constraints for Huber's theorem to be applicable. Unexpectedly, such metrics are more rigid than we initially anticipated. For instance, we found that the volume density at infinity is precisely one, and the blow-down of the metric is $\mathbb{R}^n$. Specifically, in four dimensions, we derive the $L^2$-integrability of the Ricci curvature, which directly leads to the conclusion that the Pfaffian 4-form is integrable and adheres to a Gauss-Bonnet-Chern formula.  Additionally, we demonstrate that a Gauss-Bonnet-Chern formula, previously verified by Lu and Wang under the  assumption that the second fundamental form belongs to $L^4$, remains valid for $R \in L^2$. Consequently, on an orientable 4-dimensional manifold conformal to a domain in a closed manifold, Huber's Theorem holds when $R \in L^2$, if and only if the negative part of the Pfaffian 4-form is integrable.\\

\noindent{\bf MSC Classification}\hspace{1ex}  53C21, 58J05.\\

\noindent{\bf Data Availability}\hspace{1ex} Data sharing not applicable to this article as no datasets were generated or analysed during the current study.
\end{abstract}

\section{Introduction}
In the famous paper \cite{Huber}, Huber proved a remarkable result concerning the structures of
complete surfaces: every complete surface with the integrable negative part of the Gauss curvature
is conformally equivalent to a compact surface
with a finite number of points removed.
Regrettably,  this result does not extend straightforwardly to higher dimensions. For instance,  the manifold $\mathbb{T}^2\times\R$ is flat but not conformal to any closed manifold with finite
points removed.
Therefore,  a variety of generalizations of Huber's Theorem have been established under certain supplementary curvature and other geometric assumptions, as seen in \cite{Chang-Qing-Yang2,Carron-Herzlich, Fang, Carron, Ma-Qing}, and related references.

In this paper, our focus will be on a complete manifold that conforms to a domain of a closed manifold with $\int_M|R(g)|^\frac{n}{2}dV_g<+\infty$. There are very few results available regarding Huber's Theorem in this particular direction. The only known sufficient and necessary condition is the combination of Theorem 2.1 in \cite{Carron-Herzlich} with Theorem 1.2 in \cite{ACT}, which can be summarized as follows:
\begin{thm}[Carron-Herzlich, Aldana-Carron-Tapie]
Let $\Omega$ be a domain of $(M,g_0)$, a compact Riemannian manifold of
dimension $n > 2$. Assume $\Omega$ is endowed with a complete Riemannian metric $g$
which is conformal to $g_0$. Then $M\setminus\Omega$ is a finite set if and only if $\vol(B^g_r(x_0),g)=O(r^n)$ for some point $x_0$ in $\Omega$.
\end{thm}

The primary objective of this paper is to identify more geometric constraints for Huber's Theorem.
We will investigate the geometric characteristics of a metric $g$ defined on the punctured $n$-dimensional closed ball $\overline{B}\setminus{0}$, which conforms to a smooth metric $g_0$ (defined on $\overline{B}$) with finite $\|R(g)\|_{L^\frac{n}{2}}$. Contrary to our expectations, such a metric exhibits a considerably higher degree of rigidity than previously anticipated. We demonstrate that the volume density of $g$ at infinity equals  1, and the manifold blows down to an $n$-dimensional space.
Specifically, we can state the following:

\begin{thm}\label{main} Let $g_0$ be a smooth metric defined on the closed unit $n$-ball $\overline{B}$, with $n\geq 3$. Let $g=u^\frac{4}{n-2}g_0$
be a conformal metric on $\overline{B}\setminus\{0\}$.
Assuming $\|R(g)\|_{L^\frac{n}{2}(B,g)}<+\infty$ and $\vol(\overline{B}\setminus\{0\},g)=\infty$. Then, as $r\rightarrow+\infty$, the volume ratio
$$
\frac{\vol(B_r^g(x),g)}{V_{n}r^n}\rightarrow 1,
$$
and $(\overline{B}\setminus\{0\},\frac{g}{r^2},x)$ converges to $(\R^n,0)$ in the Gromov-Hausdorff
distance, where $x$ is a fixed point in $B$ and $V_n$ represents the volume of the unit Euclidean ball. Additionally, $G^{-1}u\in W^{2,p}(B_\frac{1}{2})$
for any $p<\frac{n}{2}$, where $G$ is the Green function defined by
$$
-\Delta_{g_0}G=\delta_0,\s G|_{\partial B}=0.
$$
\end{thm}

\begin{rem}
Based on Proposition \ref{finite.volume} and Corollary \ref{complete.metric}, the assumption that $\vol(\overline{B}\setminus{0},g)=\infty$ is equivalent to the completeness of $(\overline{B}\setminus{0},g)$, when $\|R(g)\|_{L^\frac{n}{2}(B,g)}<+\infty$ is satisfied.
\end{rem}

The theorem above includes lots of unexpected pieces of information. Firstly, it implies that when $(M,g)$ is conformal to a domain of a closed manifold, Huber's result holds if and only if the volume density at infinity equals the number of ends.  In addition, each end of such a manifold
is asymptotically euclidean,  which is a useful property for  verifying  whether a manifold has a finite point conformal compactification. For instance,
by setting $g=u^\frac{4}{n-2}g_{\text{euc}}$ on $M=\mathbb{R}^n\setminus \mathbb{R}^{n-k}$, where $u=\left(\sum_{i=1}^k (x^i)^2\right)^\frac{2-k}{2}$, we find $R(g)=0$, and $(M,g)$ remains noncompact for $n>k>\frac{n}{2}+1$. Nonetheless, $(M,g)$ does not satisfy Huber's Theorem, since its blow-down  is not $\R^n$. This contrasts with cases where the total Q-curvature is finite \cite{Buzano-Nguyen,Chang-Qing-Yang1,Wang}. For instance, when $g_0$ is the Euclidean metric and $u=r^{\alpha}$, the $Q$-curvature of $g$ is 0, yet the volume density at infinity can vary widely.
Secondly, it appears that when $\|R\|_{L^\frac{n}{2}}<+\infty$, conforming to a domain of a closed manifold is a quite strong assumption. For example, if we further assume $Ric(g)\geq 0$, such a manifold must be $\R^n$ (see Corollary \ref{rigidity}).

A plausible intuitive explanation for these unusual phenomena is as follows:
Firstly, we can find sequences $r_k\rightarrow 0$ and $c_k$ such that $c_kr_k^\frac{n-2}{2}u(r_kx)$ converges weakly in $W^{2,p}(\R^n\setminus{0})$ to a positive function $u'$ for any $p<\frac{n}{2}$. The function $u'$ is harmonic since the limit metric $g_\infty={u'}^\frac{4}{n-2}g_{euc}$ is scalar flat. Moreover, $(\R^n\setminus{0},g_\infty)$ should be extendable to a cone since it can be seen as a blow-down of $(B\setminus{0},g)$. However, a positive harmonic function on $\R^n\setminus{0}$ must be in the form of $a+br^{2-n}$. When $a$ and $b$ are both non-zero, $(\R^n\setminus{0},(a+br^{2-n})^\frac{4}{n-2}g_{euc})$ becomes a complete manifold with 2 ends, and is not a cone. When $b=0$, $(\R^n\setminus{0},a^\frac{4}{n-2}g_{euc})$ is not complete near $0$. Therefore, we conclude that $u'=b|x|^{2-n}$ with $b>0$, which implies that $g_\infty$ is a flat metric defined on $\R^n\setminus{0}$.
\\

Theorem \ref{main} has a number of interesting  corollaries. First, we  examine 
a conformal map from $(\overline{B}\setminus\{0\},g_0)$ into $\R^{n+k}$. We show that if the second fundamental form
$A$ is in $L^n$ and 
the image is noncompact, then the mapping near the origin closely resembles  $\frac{x}{|x|^2}$, and the intrinsic distance is asymptotically equivalent to the distance in $\R^{n+k}$:

\begin{thm}\label{as.map}
Let $(\overline{B},g_0)$ be as in Theorem \ref{main}. Let $F:(\overline{B}\setminus\{0\},g_0)\rightarrow \R^{n+k}$ be a conformal immersion with finite $\|A\|_{L^n}$.
Suppose the volume is infinite. Then after changing the coordinates of $\R^{n+k}$, for any $r_k\rightarrow 0$ and $x_0\in B$, there exists $\lambda_k\in\R$ and $y_0\in\R^{n+k}$,
such that a subsequence of
\begin{equation}\label{MS1}
\lambda_k(F(r_kx)-F(r_kx_0))+y_0
\end{equation}
converges weakly in $W^{3,p}_{loc}(\R^n\setminus\{0\})$ to $F_\infty(x)=(\frac{x}{|x|^2},0)$ for any $p<\frac{n}{2}$. Consequently,
\begin{equation}\label{MS2}
\lim_{x\rightarrow 0}\frac{|F(x)-F(x_0)|}{d_{g_F}(x_0,x)}=1,
\end{equation}
where $g_F$ is the induced metric.
\end{thm}
Note that, if we use the coordinates change: $x\rightarrow\frac{x}{|x|^2}$,
the limit of $\lambda_k(F(r_kx)-F(r_kx_0))+y_0$
is  simply the identity map  of $\R^n$ under the new coordinates.
Therefore, Theorem \ref{as.map} can be viewed as a higher-dimensional extension of a result by S. M\"uller and V. \v{S}ver\'ak \cite[Corollary 4.2.5]{Muller-Sverak}, except that
$F$ does not have branches in our case. \\

Next, we will prove some Gauss-Bonnet-Chern formulas in 4-dimensional cases. Since the asymptotic behavior of $(\overline{B}\setminus\{0\},g)$ at
infinity is clear and simple, we can get the exact values of the error terms.

First, we discuss the formula for $Q$-curvature. For the $Q$-curvature, we use the definition in \cite{Chang1}.
Since $R(g)\in L^2$ is not strong enough to ensure the integrability of $Q$-curvature (see Example \ref{Eg3}), our first formula is stated as follows:
\begin{thm}\label{GBC1}
Let $(M_0,g_0)$ be a compact 4-dimensional orientable manifold without boundary and let $(M,g)$ be conformally equivalent to $(M_0,g_0)$ with a finite number of points removed. We assume $(M,g)$ is complete and $R(g)\in L^2(M,g)$.
Then
there exist domains $\Omega_1\subset \subset\Omega_2\subset\subset\Omega_3\cdots$, such that
$$
\bigcup_{k=1}^\infty \Omega_k=M,
$$
and
$$
\lim_{k\rightarrow+\infty}\int_{\Omega_k}Q(g)dV_g=\int_{M_0}Q(g_0)dV_{g_0}-8\pi^2m,
$$
where $m$ is the number of the ends.
\end{thm}

Under the assumptions in the above theorem, it is evident that the integrability of $Q^-$ (or $Q^+$) implies the integrability of $Q$. When $Q$ is integrable, the theorem above can be viewed as an intrinsic version of Theorem 1.5 in \cite{Lu-Wang}, where $\|A\|_{L^4}<+\infty$ was assumed. It may seem a bit unusual at first glance that such a formula is solely concerned with intrinsic properties.

As an application of the above theorem, we obtain the following:
\begin{cor}\label{R2Ric2}
Let $(M,g)$ be as in  Theorem \ref{GBC1}. Then $Riem(g)\in L^2$, where $Riem(g)$
is the curvature tensor.
\end{cor}

We denote by $Pf(g)$ the standard Pfaffian of the Riemannian metric $g$. For a closed 4-dimensional manifold  $(M_0,g_0)$, the Chern-Gauss-Bonnet formula can be expressed as
$$
\int_{M_0}Pf(g_0)=4\pi^2\chi(M_0).
$$
where $\chi(M_0)$ is the Euler characteristic of $M_0$.
It is well-known that
\begin{equation}\label{pf.int}
Pf(g)=\left(\frac{1}{8}|W(g)|^2+\frac{1}{12}R^2(g)-\frac{1}{4}|Ric(g)|^2\right)dV_{g},
\end{equation}
where $W$ is the Weyl tensor. Then
the integrability of the Pfaffian form is deduced from $L^2$-integrability of Ricci curvature and scalar curvature, along with the conformal invariance of the Weyl tensor. Furthermore, we obtain the following result:
\begin{thm}\label{GBC2}
Let $(M,g)$ and $m$ be as in  Theorem \ref{GBC1}. Then the Pfaffian of the curvature
is integrable, and
$$
\int_MPf(g)=4\pi^2\chi(M_0)-8m\pi^2,
$$
or equivalently
$$
\int_MPf(g)=4\pi^2\chi(M)-4m\pi^2.
$$
\end{thm}

We set $Pf(g)=\Phi dV_g$, where
$dV_g$ is the volume form of $g$, and define $\Phi^-dV_g$ to be the
negative part of $Pf(g)$. From the equation \eqref{pf.int},   we deduce that
$Ric(g)\in L^2(M,g)$ whenever $\Phi^-$ is integrable.
Together with Theorem 1.4 in \cite{Chen-Li}, we can establish the following

\begin{thm}\label{main2}
Let $(M,g_0)$ be a 4-dimensional oriented compact Riemannian manifold without boundary and let
$\Omega$ be a domain of $M$. Assume $\Omega$ is endowed with a complete Riemannian metric $g$
which is conformal to $g_0$ with $R(g)\in L^2(\Omega,g)$. Then $M\setminus\Omega$
is a finite set if and only if the negative part of $Pf(g)$
is integrable.
\end{thm}

This paper is organized as follows. Section 2 reviews some regularity
results of the scalar curvature equation and establishes the 3-circle Theorem. In section 3
we establish the asymptotic behaviors of the metric at infinity. Then, we  prove Theorem \ref{as.map} and Theorem \ref{GBC1}, \ref{GBC2}, \ref{main2}.
in sections 4 and 5 respectively. In the last section, we provide several examples of complete metric on the 4-dimensional punctured ball with $R\in L^2$.\\

{\bf  Acknowledgements} 
 We are grateful to the referee for  helpful corrections and highly constructive suggestions which very substantially improved 
the exposition.

\section{Preliminaries}

First, we introduce some notations that will be used throughout the remainder of the paper. We always assume $n\geq 3$ and denote
by $(B,x^1,x^2,\cdots,x^n)$ the n-dimensional unit ball, and by $B_r$ the $n$-dimensional ball of radius $r$ centered at $0$ in $\R^n$. We assume
$g_0$ is a smooth metric defined on $\overline{B}$. For simplicity,
we always assume $x^1,\cdots,x^n$ are normal coordinates of $g_0$ at $0$, then we have
\begin{equation}\label{principle.coordinates}
d_{g_0}(0,x)=|x|,\s and\s |g_{0,ij}-\delta_{ij}|\leq c|x|^2.
\end{equation}

\subsection{Regularity}

In this section, we let $v$ be a weak solution of
\begin{equation}\label{equation.epsilon}
-div(a^{ij}v_{j})=fv,
\end{equation}
where
\begin{equation}\label{aij}
0<\lambda_1\leq a^{ij},\s \|a^{ij}\|_{C^0(B_2)}+\|\nabla a^{ij}\|_{C^0(B_2)}
<\lambda_2.
\end{equation}

We have the following:
\begin{lem}\label{Lalpha}
Suppose that $v\in W^{1,2}(B_2)$ is positive and satisfies \eqref{equation.epsilon} and \eqref{aij}. We assume
$$
\int_{B_2}|f|^\frac{n}{2}\leq \Lambda.
$$
Then
$$
r^{2-n}\int_{B_r(x)}|\nabla\log v|^2<C,\s \forall B_r(x)\subset B.
$$
\end{lem}

\begin{lem}\label{regularity}
Suppose $v\in W^{1,2}(B_2)$ is positive and satisfies \eqref{equation.epsilon} and \eqref{aij}. 
Then for any $q\in (0,\frac{n}{2})$, there exists $\epsilon_0
=\epsilon_0(q,\lambda_1,\lambda_2)>0$, such that
if
$$
\int_{B_2}|f|^\frac{n}{2}<\epsilon_0,
$$
then
$$
\|\nabla\log v\|_{W^{1,q}(B)}\leq C(\lambda_1,\lambda_2,\epsilon_0),
$$
and
$$
e^{-\frac{1}{|B|}\int_{B}\log v}\|v\|_{W^{2,q}(B)}+e^{\frac{1}{|B|}\int_{B}\log v}\|v^{-1}\|_{W^{2,q}(B)}
\leq C(\lambda_1,\lambda_2,\epsilon_0).
$$
\end{lem}

For the proofs of the above two lemmas, one can refer to \cite{Dong-Li-Xu}.

\begin{cor}\label{Lp.estimate}
Suppose $v\in W^{1,2}(B_2)$ is positive and satisfies \eqref{equation.epsilon} and \eqref{aij}. Assume 
$$
\int_Bv^2<\Lambda_1.
$$ 
Then for any $q<\frac{n}{2}$ and $p>2$, there exists $C_1$ and $C_2$, such that
$$
\|v\|_{W^{2,q}(B)}<C_1,\s and\s \int_Bv^p\leq \|v\|_{L^2(B)}^2+C_2\|v\|_{L^2(B)}.
$$ 
\end{cor}

\proof
Note that for any $E\subset B$,
\begin{equation}\label{log+}
\int_E(\log v)^+ = \int_{E\cap\{v>1\}}\log v
\leq\int_E v^2 < C.
\end{equation}
Utilizing \eqref{regularity}, we obtain that $
\|v\|_{W^{2,q}(B)}<C.$

Select $q<\frac{n}{2}$, such that $W^{2,q}(B)$ can be embedded into $L^{2p}(B)$.
Since
$$
|\{v\geq 1\}\cap B|\leq \int_{B}v^2,
$$
we have
\begin{eqnarray*}
\int_{B}v^{p}&\leq&\int_{\{v\leq 1\}\cap B}v^2+
\left(\int_{\{v\geq 1\}\cap B}v^{2p}\right)^\frac{1}{2}|\{v\geq 1\}\cap B|^\frac{1}{2}\\
&\leq& \int_{B}v^2+\left(\int_{B}v^{2p}\right)^\frac{1}{2}\left(\int_{B}v^2\right)^\frac{1}{2}\leq \int_{B}v^2+C\|v\|_{W^{2,q}(B)}^p\left(\int_{B}v^2\right)^\frac{1}{2}.
\end{eqnarray*}
\endproof

\subsection{Convergence of Distance Functions} The distance between two points $x$ and $y$ on a manifold $(M,g)$ is defined
as the infimum of the lengths of piecewise smooth
curves joining them.
We will use the following proposition:
\begin{pro}\label{distance}
Let $a_k\rightarrow 0^+$ and $g_k=g_{k,ij}dx^i\otimes dx^j$ be a smooth metric defined on $B_\frac{1}{a_k}\setminus B_{a_k}$.
Assume $g_k$ and $g_k^{-1}$ converge to $g_{euc}$ in $W^{1,p}_{loc}(\R^n\setminus\{0\})$
for any $p\in(n-1,n)$. Then after passing to a subsequence,
$d_{g_k}$ converges to $d_{g_{euc}}$ in $C^0((B_r\setminus B_\frac{1}{r})\times (B_r\setminus B_\frac{1}{r}))$
for any $r>1$.
\end{pro}

\proof The arguments in \cite[Section 3]{Dong-Li-Xu} use properties of complete metrics, so they
can not  be applied here directly. For this reason, we let $t> 1$ and take nonnegative $\phi_t\in C^\infty(\R)$, which satisfies:
1). $\phi_t$ is 1 on $[\frac{1}{t},t]$ and 0
on $(-\infty,\frac{1}{2t}]\cup[2t,+\infty)$; 2). $|\phi'|<4t$.
Define
$$
\hat{g}_{k,t}=\phi_t(|x|) g_k+(1-\phi_t(|x|))g_{euc}.
$$
Obviously, $\hat{g}_{k,t}$ is complete on $\R^n$.

We have
\begin{eqnarray*}
\int_{B_{2t}\setminus B_{t}}|\nabla({\hat g}_{k,t}-g_{euc})|^pdx&\leq&
Ct^p\int_{B_{2t}\setminus B_{t}}|g_{k}-g_{euc}|^pdx+C\int_{B_{2t}\setminus B_{t}}|\nabla(g_{k}-g_{euc})|^pdx\\
&\leq&C(t)\|g_{k}-g_{euc}\|^p_{W^{1,p}(B_{2t}\setminus B_{t})}.
\end{eqnarray*}
A similar estimate can be obtained on $B_\frac{1}{t}\setminus B_\frac{1}{2t}$ using the same argument.

Note that
$$
det(\phi_tg_k+(1-\phi_t)g_{euc})\geq \prod_{i=1}^n(\phi_t a_i+(1-\phi_t)),
$$
where $a_1$, $\cdots$, $a_n$ are eigenvalues of $g_k$. Then
$$
det(\phi_tg_k+(1-\phi_t)g_{euc})\geq \frac{1}{2^n}\min\{det(g_k),1\},
$$
which implies that
$$
\frac{1}{det(\hat{g}_{k,t})}\leq C(1+\frac{1}{det(g_k)})=C(1+det(g_k^{-1})).
$$
Since $g_k^{-1}$ converges in $W^{1,p}_{loc}(\R^n\setminus\{0\})$ for
any $p\in(n-1,n)$, $(g_{k}^{-1})$ is bounded in $L^q(B_r\setminus B_\frac{1}{r})$ for any $q$. Then $1/det(\hat{g}_{k,t})$ is bounded in $L^q(B_r)$ for any $q$. Recall that
the inverse of a matrix is just the adjugate matrix divided by the determinant. Then $\hat{g}^{-1}_{k,t}$
is bounded in $W^{1,p}(B_r)$.

Note that $\hat{g}_{k,t}=g_{ecu}$ on $B_{2t}^c$. It is not difficult to check that for any fixed $r$, there exists $r'$, such that
any geodesic between two points $x$, $y\in B_r$ must lie in $B_{r'}$.
Then, using the arguments in \cite[Section 3]{Dong-Li-Xu}, a subsequence of $d_{\hat{g}_{k,t}}$ converges to $d_{g_{euc}}$
in $C^0(B_r\times B_r)$ for any $r$.
Thus, after passing to a subsequence, we can find $t_k\rightarrow+\infty$, such that $d_{\hat{g}_{k,t_k}}$
converges to $d_{g_{euc}}$ in $C^0(B_r\times B_r)$ for any $r$. For simplicity, we set $\tilde{g}_k=\hat{g}_{k,t_k}$ and assume $d_{\tilde{g}_k}$ converges in $C^0(B_r\times B_r)$ for any $r$.\\

Now, we start to prove that $d_{g_k}$ converges to $d_{g_{euc}}$ in $C^0((B_r\setminus B_\frac{1}{r})\times (B_r\setminus B_\frac{1}{r}))$.

Let $\lambda_k(x)$ be the lowest eigenvalue of $g^{ij}_k$.
Since $|\nabla^{g_k}_xd_{g_k}(x,y)|\leq 1$ for a.e. $x$, we have
$$
|\nabla_xd_{g_k}(x,y)|^2\leq\frac{1}{\lambda_k(x)} \leq C\sum_{ij}|g_{k,ij}(x)|,
$$
which implies that $d_{g_k}$ is bounded $W^{1,q}((B_r\setminus B_\frac{1}{r})\times (B_r\setminus B_\frac{1}{r}))$ for any $r$ and $q>0$. Then, we may assume $d_{g_k}$ converges to a function $d$ in
$C^0((B_r\setminus B_\frac{1}{r})\times (B_r\setminus B_\frac{1}{r}))$ for any $r$. By the trace embedding theorem, for any $x$, $y\in B_r\setminus B_\frac{1}{r}$, we have
$$
d(x,y)\leq d_{g_{euc}}(x,y).
$$

Next, we show $d(x,y)\geq d_{g_{euc}}(x,y)$.
Let $\gamma_k$ be a curve from $x$ to $y$ in $\R^n\setminus\{0\}$, such that
$$
L_{g_k}(\gamma_k)\leq d_{g_k}(x,y)+\frac{1}{k}.
$$

Let $\lambda>d(x,y)+r+1$. We claim that $\gamma_k\subset B_{\lambda}$ when
$k$ is sufficiently large. Suppose that $\gamma_k\cap \partial B_{\lambda}\neq\emptyset$.
It is easy to check that
$$
d_{g_k}(\partial B_\lambda,\partial B_r)\leq
L_{g_k}(\gamma_k)\leq d_{g_k}(x,y)+\frac{1}{k}\rightarrow d(x,y).
$$
However,
$d_{g_k}(\partial B_\lambda,\partial B_r)=d_{\tilde g_k}(\partial B_\lambda,\partial B_r)$ when $k$
is sufficiently large, and
$$
\lim_{k\rightarrow+\infty}d_{\tilde g_k}(\partial B_\lambda,\partial B_r)=
d_{ g_{euc}}(\partial B_\lambda,\partial B_r)=\lambda-r,
$$
which leads to a contradiction.

The rest of the proof can be divided into 2 cases. Case 1, we assume $\gamma_k\cap\partial B_\frac{1}{t_k}=\emptyset$.
In this case,
$$
L_{g_k}(\gamma_k)=L_{\tilde{g}_k}(\gamma_k)\geq d_{\tilde{g}_k}(x,y).
$$
Case 2, we may assume $x_k$ and $y_k$ to be the first and the last point in $\gamma_k\cap \partial B_\frac{1}{t_k}$ respectively. Then
\begin{eqnarray*}
L_{g_k}(\gamma_k)
&\geq& L_{g_k}(\gamma_k|_{[x,x_k]})+ L_{g_k}(\gamma_k|_{[y_k,y]})\\
&=&L_{\tilde{g}_k}(\gamma_k|_{[x,x_k]})+ L_{\tilde{g}_k}(\gamma_k|_{[y_k,y]})\\
&\geq&d_{\tilde{g}_k}(x,x_k) +d_{\tilde{g}_k}(y_k,y)\\
&\geq& d_{\tilde{g}_k}(x,y)-
d_{\tilde{g}_k}(x_k,y_k).
\end{eqnarray*}
Thus, for both cases, we have
$$
d(x,y)=\lim_{k\rightarrow+\infty}d_{g_k}(x,y)\geq \lim_{k\rightarrow+\infty}d_{\tilde{g}_k}(x,y)=d_{g_{euc}}(x,y).
$$
\endproof

\subsection{Three Circles Theorem}
In this section, we present the Three Circles Theorem. It is convenient to state and prove this theorem
on pipes. We let $Q=[0,3L]\times S^{n-1}$, and
$$
Q_i=[(i-1)L,iL]\times S^{n-1},\s i=1,2,3.
$$
Set $g_Q=dt^2+g_{S^{n-1}}$ and $dV_Q=
dV_{g_Q}$.

We first state this theorem for the case of $g=g_Q$ and $R=(n-1)(n-2)$:
\begin{lem}\label{3-circle.R=0}
Let $u\neq 0$ solve the following equation on $Q$:
$$
-\Delta u+\frac{(n-2)^{2}}{4}u=0.
$$
Then there exists $L_0$, such that for any $L>L_0$, we have

\begin{itemize}
\item[1) ]$\int_{Q_1}u^2dV_{Q}\leq e^{-L}
\int_{Q_2}u^2dV_{Q}$ implies
$\int_{Q_2}u^2dV_{Q}< e^{-L}\int_{Q_3}u^2dV_{Q};$

\item[2) ] $\int_{Q_3}u^2dV_{Q}\leq e^{-L}
\int_{Q_2}u^2dV_{Q}$ implies
$\int_{Q_2}u^2dV_{Q}< e^{-L}\int_{Q_1}u^2dV_{Q};$

\item[3) ]either
$\int_{Q_2}u^2dV_{Q} < e^{-L}\int_{Q_1}u^2dV_{Q}$
or
$\int_{Q_2}u^2dV_{Q}< e^{-L}\int_{Q_3}u^2dV_{Q}.
$
\end{itemize}
\end{lem}

For the proof, one can refer to \cite{Bru,Li-Zhou}.
Next, we discuss the general case.
\begin{thm}\label{3-circle}
Let $g_0$ be a metric over $Q$ and $u\in W^{2,p}$
which solves the equation
$$
-\Delta_{g_0}u+c(n)R(g_0)u=fu.
$$
Then for any $L>L_0$, there exist $\epsilon_0'$, $\tau$, such that
if
\begin{equation}\label{3circle.assumption}
\|g_0-g_Q\|_{C^2(Q)}<\tau,\s \int_Q|f|^\frac{n}{2}dV_{g_0}<\epsilon_0',
\end{equation}
then

\begin{itemize}
\item[1) ]$\int_{Q_1}u^2dV_{g_0}\leq e^{-L}
\int_{Q_2}u^2dV_{g_0}$ implies
$\int_{Q_2}u^2dV_{g_0}\leq e^{-L}\int_{Q_3}u^2dV_{g_0};$

\item[2) ] $\int_{Q_3}u^2dV_{g_0}\leq e^{-L}
\dint_{Q_2}u^2dV_{g_0}$ implies
$\int_{Q_2}u^2dV_{g_0}\leq e^{-L}\int_{Q_1}u^2dV_{g_0};$

\item[3) ]either
$
\int_{Q_2}u^2dV_{g_0} \leq e^{-L}\int_{Q_1}u^2dV_{g_0}$
or $\int_{Q_2}u^2dV_{g_0}\leq e^{-L}\int_{Q_3}u^2dV_{g_0}.
$
\end{itemize}
\end{thm}

\proof
If the statement in 1) is false for an $L>L_0$, we can find $g_k$,
$u_k$ and $f_k$, s.t.
$$
g_k\rightarrow g_Q \mbox{ in }C^2(Q),\s \int_{Q}|f_k|^\frac{n}{2}dV_{g_k}
\rightarrow 0,$$
$$
-\Delta_{g_k}u_k+c(n)R(g_k)u_k=f_ku_{k},
$$
and
$$
\int_{Q_1}u_k^2dV_{g_k}\leq e^{-L}\int_{Q_2}u_k^2dV_{g_k},\s
\int_{Q_2}u_k^2dV_{g_k} > e^{-L}\int_{Q_3}
u_k^2dV_{g_k}.
$$
Let $$v_k=\frac{u_k}{\|u_k\|_{L^2(Q_2,g_k)}}.$$
We have
$$
\int_{Q_1}v_k^2dV_{g_k}\leq e^{-L}\int_{Q_2}v_k^2dV_{g_k},\s
\int_{Q_3}v_k^2dV_{g_k} < e^{L}\int_{Q_2}v_k^2dV_{g_k},
$$
and
$$
\int_{Q_2}v_k^2dV_{g_k}=1.
$$
Thus
$$
\int_Qv_k^2dV_{g_k}\leq
C.
$$
$v_k$ satisfies
$$
-\Delta_{g_k} v_k+c(n)R(g_k)v_k=f_kv_k.
$$
By Corollary \ref{Lp.estimate} and Sobolev embedding theorem, $v_k$ converges to a function $v$ in $W^{1,2}_{loc}$, where $v$ satisfies:
$$
-\Delta v+\frac{(n-2)^2}{4}v=0,\s and\s
\int_{Q_2}|v|^2dV_Q=1.
$$
Thus $v\neq 0$.

Moreover,
$$
\int_{[\epsilon,L]\times S^{n-1}}v^2dV_Q=\lim_{k\rightarrow+\infty}
\int_{[\epsilon,L]\times S^{n-1}}v_k^2dV_{g_k}\leq e^{-L}
\lim_{k\rightarrow+\infty}\int_{Q_2}v_k^2dV_{g_k},
$$
letting $\epsilon \rightarrow 0$ gives
$$
\int_{Q_1}v^2dV_Q\leq e^{-L}\int_{Q_2}v^2dV_Q.
$$
Similarly, there holds
$$
\int_{Q_3}v^2dV_Q\leq e^{L}\int_{Q_2}v^2dV_Q,
$$
which contradicts Lemma \ref{3-circle.R=0}. Hence, the statements in 1) are proved. Using the same arguments, we can easily carry out the proof of 2) and 3).
~\endproof

\begin{thm}\label{3circle}
Let $g=u^\frac{4}{n-4}g_0$ be a smooth metric
defined on $\overline{B}\setminus\{0\}$ with
$$
\int_{B}|R(g)|^\frac{n}{2}dV_g<+\infty,\s \vol(\overline{B}\setminus\{0\},g)=+\infty.
$$
Then for any $\vartheta>e^{L_0}$ there exists $r_0$, such that
for any $r<r_0$, there holds
\begin{equation}\label{3.circle.L2}
\int_{B_{r}\setminus B_{r\vartheta^{-1}}}\frac{u^2}{|x|^2}dV_{g_0}\leq \frac{1}{\vartheta}
\int_{B_{r\vartheta^{-1}}\setminus B_{r\vartheta^{-2}}}\frac{u^2}{|x|^2}dV_{g_0}.
\end{equation}
Moreover, we have
\begin{equation}\label{infinite}
\lim_{k\rightarrow+\infty}\int_{B_{\vartheta^{-k}r}\setminus B_{\vartheta^{-k-1}r}}\frac{u^2}{|x|^2}dV_{g_0}=
+\infty\s and\s\lim_{k\rightarrow+\infty}\int_{B_{\vartheta^{-k}r}\setminus B_{\vartheta^{-k-1}r}}u^\frac{2n}{n-2}dV_{g_0}=+\infty.
\end{equation}
\end{thm}

\proof
Put
\begin{align*}
\phi(t,\theta)=(e^{-t},\theta),
\end{align*}
and
\begin{align*}
g'(t,\theta)=\phi^*(g)=v^\frac{4}{n-2}\hat{g}(t,\theta),
\end{align*}
where $\hat{g}(t,\theta)=e^{2t}\phi^*(g_0)$, which converges to $dt^2+g_{{S}^{n-1}}$ as $t\rightarrow+\infty$. Then
\begin{align*}
v^\frac{4}{n-2}(t,\theta)=u^\frac{4}{n-2}(e^{-t},\theta)e^{-2t},
\end{align*}
\begin{align*}
-\Delta_{\hat{g}}v+c(n)R(\hat{g})v=c(n)R(g')v^\frac{n+2}{n-2}:=fv,
\end{align*}
\begin{align*}
\int_{[a,b]\times S^{n-1} }|f|^\frac{n}{2}dV_{\hat{g}}=c\int_{B_{e^{-a}}\setminus B_{e^{-b}}}|R({g})|^\frac{n}{2}dV_g,
\end{align*}
and
$$
\int_{B_r\setminus B_{r/\vartheta}}\frac{u^2}{|x|^2}dV_{g_0}=\int_{[-\log r,-\log r+\log \vartheta]\times S^{n-1}}v^2dV_{\hat{g}}.
$$

Without loss of generality, we assume 
\begin{equation}\label{small.assumption}
\|\hat g-g_Q\|_{C^2([0,+\infty)\times S^{n-1})}<\tau,\s \int_{[0,+\infty)\times S^{n-1}}|R|^\frac{n}{2}dV_{g}<\epsilon_0'.
\end{equation}

Suppose \eqref{3.circle.L2} is not true, i.e., we can find $r_k\rightarrow 0$, such that 
\begin{equation}\label{vol.infinity.assu}
\int_{B_{r_k\vartheta^2}\setminus B_{r_k\vartheta}}\frac{u^2}{|x|^2}dV_{g_0}> \frac{1}{\vartheta}
\int_{B_{r_k\vartheta }\setminus B_{r_k}}\frac{u^2}{|x|^2}dV_{g_0}.
\end{equation}
We set
$$
\Omega_{k,m}=[-\log r_k-(m_k-m+1)\log \vartheta,-\log r_k-(m_k-m)\log \vartheta]\times S^{n-1},
$$
where $m=1$, $\cdots$, $m_k=\left[\frac{-\log r_k}{\log\vartheta}\right]$.
Then \eqref{vol.infinity.assu} is equivalent to
$$
\int_{\Omega_{k,m_k}}v^2dV_{\hat{g}}\leq \frac{1}{\vartheta}\int_{\Omega_{k,m_k-1}}v^2dV_{\hat{g}}.
$$
By Theorem \ref{3-circle}, we get
$$
\int_{\Omega_{k,m_k-1}}v^2dV_{\hat{g}}\leq \frac{1}{\vartheta}\int_{\Omega_{k,m_k-2}}v^2dV_{\hat{g}}.
$$
Step by step, we get
$$
\int_{\Omega_{k,m}}v^2dV_{\hat{g}}\leq \vartheta^{-(m-1)}\int_{\Omega_{k,1}}v^2dV_{\hat{g}}\leq \vartheta^{-(m-1)}\int_{[0,2\log\vartheta]\times S^{n-1}}v^2dV_{\hat{g}}\leq C\vartheta^{-m}.
$$
By Corollary \ref{Lp.estimate},
\begin{eqnarray*}
\int_{\Omega_{k,m}}v^\frac{2n}{n-2}dV_{\hat{g}}
&\leq& C(\vartheta^{-m}+\vartheta^{-\frac{m}{2}}),
\end{eqnarray*}
hence
\begin{eqnarray*}
\int_{B\setminus B_{r_k}}u^\frac{2n}{n-2}dV_{g_0}&\leq& C\int_{B\setminus B_{r_k}}u^\frac{2n}{n-2}dx\leq C\int_{B\setminus B_{\vartheta^{-1}}}u^\frac{2n}{n-2}dx+ C\sum_{m=1}^{m_k}\int_{B_{r_k\vartheta^{m}}\setminus B_{r_k\vartheta^{m-1}}}u^\frac{2n}{n-2}dx\\
&\leq& \sum_m\frac{C}{\vartheta^\frac{m}{2}}<C(\vartheta),
\end{eqnarray*}
where $C(\vartheta)$ is independent of $k$. Letting $k\rightarrow \infty$, we get a contradiction.

Thus, we get \eqref{3.circle.L2},
which implies from Theorem \ref{3-circle} that
$$
\int_{B_{r\vartheta^{-m}}\setminus B_{r\vartheta^{-m-1}}}\frac{u^2}{|x|^2}dV_{g_0}\geq C\vartheta^m.
$$
Since
$$
\int_{B_{r\vartheta^{-m}}\setminus B_{r\vartheta^{-m-1}}}\frac{u^2}{|x|^2}dV_{g_0}
\leq C(\log \vartheta)^\frac{2}{n}(\int_{B_{r\vartheta^{-m}}\setminus B_{r\vartheta^{-m-1}}}u^\frac{2n}{n-2}dV_{g_0})^\frac{n-2}{n},
$$
we get
$$
\lim_{k\rightarrow +\infty}\int_{B_{r\vartheta^{-m}}\setminus B_{r\vartheta^{-m-1}}}u^\frac{2n}{n-2}dV_{g_0}=\infty.
$$
\endproof

\begin{pro}\label{finite.volume}
Let $g=u^\frac{4}{n-4}g_0$ be a smooth metric
defined on $\overline{B}\setminus\{0\}$ with
$$
\int_{B}|R(g)|^\frac{n}{2}dV_g<+\infty,\s \vol(\overline{B}\setminus\{0\},g)<+\infty.
$$
Then $(\overline{B}\setminus\{0\},g)$ is bounded.
\end{pro}

\proof
We need to show that there exist $r$ and $C$, such that for any $x$  sufficiently close to $0$, we can find $x'\in \bar B_1\setminus B_r$, such have $d_g(x,x')<C$.

Let $g'$, $\hat g$, $v$ be as in the proof of  Theorem \ref{3circle} and assume \eqref{small.assumption} holds.
Set
$$
\Omega_{m}=[-\log r_0+(m-1)\log \vartheta,-\log r_0+m\log \vartheta]\times S^{n-1},
$$
where $-\log r_0\in [0,\log\vartheta)$ such that $x=(-\log r_0+m_0\log\vartheta,\theta)$ for some $m_0\in \mathbb{Z}^+$ and $\theta\in S^{n-1}$. 

By Theorem \ref{3-circle}, if there exists $m$, such that
$\int_{\Omega_{m}}v^2dV_{\hat{g}}\leq \vartheta^{-1}\int_{\Omega_{m+1}}v^2dV_{\hat{g}}$,
then 
$$
\int_{\Omega_{m+m'}}v^2dV_{\hat{g}}\geq C\vartheta^{m'}\rightarrow+\infty,
$$ 
which is impossible, since
$$
\int_{\Omega_{m+m'}}v^2dV_{\hat{g}}\leq C\left(\int_{\Omega_{m+m'}}v^\frac{2n}{n-2}dV_{\hat{g}}\right)^\frac{n-2}{n}\leq C(\vol(B\setminus\{0\},g))^\frac{n-2}{n}
<+\infty.
$$
Then
$$
\int_{\Omega_{m+1}}v^2dV_{\hat{g}}\leq \vartheta^{-1}\int_{\Omega_{m}}v^2dV_{\hat{g}},
$$
which implies that
$$
\int_{\Omega_{m}}v^2dV_{\hat{g}}<C\vartheta^{-m}.
$$
By Corollary \ref{Lp.estimate}, $\|v\|_{L^p(\Omega_m)}<C(p)\vartheta^{-m/(2p)}$ for any $p$,
hence, for any $q<\frac{n}{2}$, $\|\Delta_{\hat g} v\|_{L^q(\Omega_m)}<C(q)\vartheta^{-a(q)m}$ for some $a(q)>0$. By the standard elliptic estimate, we get
$\|v\|_{W^{2,q}(\Omega_m)}<C\vartheta^{-a(q)m}$. It follows from the Sobolev inequality that for any $q'\in (n-1,n)$, the inequality 
$\|v^\frac{2n}{n-2}\|_{W^{1,q'}(\Omega_m)}<C(q')\vartheta^{-a(q')m}$ holds for some positive constants
$C(q')$ and $a(q')$. 

For convenience, we set 
$$
t_m=-\log r_0+m\log\vartheta,\s x_m=(t_m,\theta).
$$
By the classical trace embedding theorem (cf. \cite[Theorem 4.12]{AF}), we have
\begin{eqnarray*}
d_{\hat g}(x_m,x_{m+1})&\leq& C\int_{0}^{\log\vartheta} v^\frac{2}{n-2}(t+t_m,\theta)dt\leq C(\int_0^{\log\vartheta}v^\frac{2n}{n-2}(t+t_m,\theta))^\frac{1}{n}\\
&\leq& C\|v^\frac{2n}{n-2}\|^\frac{1}{n}_{W^{1,q'}(\Omega_m)}\leq
C\vartheta^{-\frac{a(q')m}{n}}. 
\end{eqnarray*}
Then 
$$
d(x,x_0)<\sum_{m=1}^{m_0}d(x_m,x_{m-1})<C.
$$
~\endproof

\section{asymptotic properties}

In this section, we always assume that \eqref{principle.coordinates} holds and $g=u^\frac{4}{n-4}g_0$  denotes a smooth metric
on $\overline{B}\setminus\{0\}$ with
$$
\vol(B\setminus\{0\},g)=\infty.
$$

First of all, we prove the following lemma:
\begin{lem}\label{asy}
Let $r_k\rightarrow 0$. After passing to a subsequence, we can find $c_k>0$, such that
$c_kr_k^{\frac{n-2}{2}}u(r_kx)$ converges to $|x|^{2-n}$ weakly in $W^{2,p}_{loc}(\R^n\setminus\{0\})$
for any $p\in[1,\frac{n}{2})$. Moreover, $d_{g_k}$ converges to
$d_{|x|^{-4}g_{euc}}$ in $C^0((B_r\setminus B_\frac{1}{r})\times(B_r\setminus B_\frac{1}{r}))$ for any $r>1$, where
$g_{k,ij}=r_k^2(c_ku(r_kx))^\frac{4}{n-2}g_{0,ij}(r_kx)$.
\end{lem}
\proof
Define
$$
u_k(x)=r_k^{\frac{n-2}{2}}u(r_kx)c_k,
$$
where $c_k$ is chosen such that $\int_{\partial B_1}\log u_k d{{S}^{n-1}}=0$. It is easy to check that
$-\Delta_{g_k} u_k=f_ku_k$, where
$$
f_k=-c(n)r_k^2R_{g_0}(r_kx)+c(n)R(r_kx)(r_k^\frac{n-2}{2}u(r_kx))^\frac{4}{n-2},$$
and for sufficiently large $k$,
$$
\|f_k\|_{L^\frac{n}{2}(B_r\setminus B_\frac{1}{r})}\leq Cr_k^2+C
\left(\int_{B_{rr_k}\setminus B_\frac{r_k}{r}}|R|^\frac{n}{2}u^\frac{2n}{n-2}dx\right)^\frac{2}{n}
\leq \min\{\epsilon_0,\epsilon_0'\}.
$$

By Lemma \ref{Lalpha} and the Poincar\'e inequality (c.f. \cite[Theorem 5.4.3]{Attouch-Buttazzo-Michaille}),
$\log u_k$ is bounded in $W^{1,2}(B_r\setminus B_\frac{1}{r})$.
Then $u_k$ is bounded in $W^{2,p}(B_r\setminus B_\frac{1}{r})$ by using Lemma \ref{regularity}. Thus $u_k$ converges weakly to a positive harmonic function $u'$
locally on $\R^n\setminus\{0\}$ with $\int_{\partial B_1}\log u'd{S}^{n-1}=0$.
According Corollary 3.14 in \cite{Axler-Bourdon-Ramey}, $u'$ is written as
$$
u'=a+b|x|^{2-n},
$$
where $a$ and $b$ are nonnegative constants.
Applying Theorem \ref{infinite}, we get
$$
\int_{B_{r_kr}\setminus B_{r_kr\vartheta^{-1}}}\frac{u^2}{|x|^2}dV_{g_0}\leq \frac{1}{\vartheta}
\int_{B_{r_kr\vartheta^{-1}}\setminus B_{r_kr\vartheta^{-2}}}\frac{u^2}{|x|^2}dV_{g_0},
$$
which implies that
$$
\int_{B_r\setminus B_{r\vartheta^{-1}}}\frac{u_k^2}{|x|^2}dV_{g_{0,k}}\leq \frac{1}{\vartheta}
\int_{B_{r\vartheta^{-1}}\setminus B_{r\vartheta^{-2}}}\frac{u_k^2}{|x|^2}dV_{g_{0,k}},
$$
where $g_{0,k}=g_{0,ij}(r_kx)dx^i\otimes dx^j$.
Taking the limit, we obtain
$$
\int_{B_r\setminus B_{r\vartheta^{-1}}}\frac{{u'}^2}{|x|^2}dx\leq \frac{1}{\vartheta}
\int_{B_{r\vartheta^{-1}}\setminus B_{r\vartheta^{-2}}}\frac{{u'}^2}{|x|^2}dx.
$$
Letting $r$ be sufficiently large, we get $a=0$.
Since $\int_{\partial B_1}\log u'=0$, $b=1$.

By changing coordinates: $x\rightarrow \frac{x}{|x|^2}$, we see that ${u'}^\frac{4}{n-2}g_{euc}$
is just $g_{euc}$ in the new coordinates. The convergence of $d_{g_k}$ follows
from Proposition \ref{distance} directly.
\endproof

In the preceding lemma, we did not express $c_k$ in terms of $r_k$, which limits our understanding of the behavior of $u(r_kx)$. Nonetheless, the lemma is sufficiently strong to derive the following decay properties:

\begin{cor}\label{decay}
For any $\tau\in(0,1)$, there exists $\delta$ such that for any $r<\delta$, the following hold:
\begin{eqnarray}\label{volume.decay}
(1-\tau)\frac{1}{2^n}\leq \frac{\vol(B_{2r}\setminus B_r,g)}{\vol(B_r\setminus B_{r/2},g)}
\leq \frac{1}{2^n}(1+\tau);
\\
\label{volume.over.distance}
(1-\tau)(2^n-1)\leq
\frac{\vol(B_{2r}\setminus B_r,g)}{V_n (d_g(2rx_0,rx_0))^n}
\leq (2^n-1)(1+\tau),\s \forall x_0\in\partial B;
\\
\label{nabla.decay}
\frac{\int_{B_{r}\setminus B_{r/2}}|x|^\beta|\nabla_{g_0} u|^\alpha dV_{g_0}}{\int_{B_{2r}\setminus B_{r}} |x|^\beta|\nabla_{g_0} u|^\alpha dV_{g_0}}<2^{(n-1)\alpha-n-\beta}(1+\tau),
\s \forall \alpha\in[1,n), \s \beta\in\R;\\
\label{u.decay}
\frac{\int_{B_{r}\setminus B_{r/2}}|x|^\beta u^\alpha dV_{g_0}}{\int_{B_{2r}\setminus B_{r}} |x|^\beta u^\alpha dV_{g_0}}<
2^{(n-2)\alpha-n-\beta}(1+\tau),\s \forall \alpha\in[1,+\infty),\s \beta\in\R;\\
\label{logu.decay}
\frac{1}{2^{n-2}}(1-\tau)\leq\frac{\int_{\partial B_r}|\nabla_{g_0}\log u|dS_{g_0}}{\int_{\partial B_{2r}}|\nabla_{g_0}\log u|dS_{g_0}}
\leq \frac{1}{2^{n-2}}(1+\tau);\\
\label{distance.decay}
2^{-1}(1-\tau)\leq\frac{d_g(2rx_0,rx_0)}{d_g(rx_0,rx_0/2)}\leq 2^{-1}(1+\tau),\s \forall x_0\in\partial B;\\
\label{point.circle.decay}
(1-\tau)\leq\frac{d_g(rx_0,\partial B_{2r})}{d_g(rx_0,2rx_0)}\leq (1+\tau),\s \forall x_0\in\partial B;\\
\label{sphere.diam.decay}
2^{-1}(1-\tau)\leq\frac{\diam(\partial B_{2r},g)}{\diam(\partial B_r,g)}\leq 2^{-1}(1+\tau);\\
\label{sphere.distance.decay}
2^{-1}(1-\tau)\leq\frac{d_g(\partial B_{2r},\partial B_{r})}{d_g(\partial B_{r},\partial B_{r/2})}\leq 2^{-1}(1+\tau);\\
\label{annulus.diam.decay}
4(1-\tau)\leq\frac{\diam( B_{2r}\setminus B_r,g)}{d_g(\partial B_{2r},\partial B_{r})}\leq 4(1+\tau).
\end{eqnarray}
\end{cor}

\proof
Let $u_k=c_kr_k^\frac{n-2}{2}u(r_kx)$ be as in the proof of Lemma \ref{asy}, which converges to $u'=|x|^{2-n}$
weakly in $W^{2,p}_{loc}(\R^n\setminus\{0\})$ for any $p<\frac{n}{2}$.
Then we may assume $\nabla u_k$ converges in $L^q_{loc}(\R^n\setminus\{0\})$
for any $q<n$, and $u_k$ converges in $L^q_{loc}(\R^n\setminus\{0\})$ for any
$q>0$. By the trace inequality, we can also assume
$\log u_k$ converges in $L^1(\partial B_t)$.

Now, we prove the right-hand side inequality of \eqref{volume.decay}:
assume it is not valid, then there exists $r_k\rightarrow 0$, such that
$$
\frac{\vol(B_{2r_k}\setminus B_{r_k},g)}{\vol(B_{r_k}\setminus B_{r_k/2},g)}> \frac{1}{2^n}(1+\tau),
$$
which means that
$$
\frac{\vol(B_{2}\setminus B_1,g_k)}{\vol(B_{1}\setminus B_{1/2},g_k)}> \frac{1}{2^n}(1+\tau),
$$
where $g_{k,ij}=u_k^\frac{4}{n-2}g_{0,ij}(r_kx)$. Letting $k\rightarrow+\infty$,
we get
$$
\frac{\vol(B_{2}\setminus B_1,g_\infty)}{\vol(B_{1}\setminus B_{1/2},g_\infty)}\geq \frac{1}{2^n}(1+\tau),
$$
where $g_\infty=|x|^{2-n}g_{euc}$. A contradiction.

Since the proofs of other inequalities are almost the same, we omit them.
\endproof

The inequalities \eqref{volume.decay}-\eqref{annulus.diam.decay} will be used to estimate quantities on $B_{2r}\setminus B_r$. For instance, using \eqref{volume.decay}, we have
$$
\vol(B_{2^kr}\setminus B_{2^{k-1}r},g)\leq \left(\frac{1+\tau}{2^n}\right)^{k-1}\vol(B_{2r}\setminus B_r,g),
$$
which implies 
\begin{equation}\label{volume.estimate}
\vol(B\setminus B_r,g)\leq C\vol(B_{2r}\setminus B_{r},g).
\end{equation}
We provide  several additional applications.

\begin{cor}\label{complete.metric}
The manifold $(\overline{B}\setminus\{0\},g)$ is complete.
\end{cor}
\proof
Consider a sequence $\{x_k\}$ that does not contain a convergent subsequence. Then
$x_k$ converges to $0$. Let $x_0\in \partial B$ be fixed. To show completeness, it suffices to prove that
$d_g(x_k,x_0)\rightarrow+\infty$.

By \eqref{distance.decay} and \eqref{point.circle.decay}, we have $d_g(x_k,2x_k)\rightarrow+\infty$,
and
$$
d_g(x_0,x_k)\geq d_g(x_k,\partial B_{2|x_k|})\geq(1-\tau)d_g(x_k,2x_k)\rightarrow+\infty.
$$
\endproof

\begin{cor}
Let $r_k$, $c_k$ be as in Lemma \ref{asy}.
Let $x'\in\partial B$ and $\rho_k=d_g(x',r_kx')$. Then
$$
\lim_{k\rightarrow+\infty}c_k^\frac{2}{n-2}\rho_k=1.
$$
\end{cor}
\proof
Let $u_k$ and $g_k$ be as in Lemma \ref{asy}. Set $g_\infty=|x|^{-4}g_{euc}$.
By Lemma \ref{asy}, for any $\sigma=2^m$, we have
$$
d_{g_k}(x',\sigma x')\rightarrow
d_{g_\infty}(x',\sigma x')=1-1/\sigma.
$$
Thus
$$
c_k^\frac{2}{n-2}d_{g}(r_kx',r_k\sigma x')=d_{g_k}(x',\sigma x')\rightarrow 1-1/\sigma.
$$
Using \eqref{distance.decay}, we get
$$
c_k^\frac{2}{n-2}d_g(r_k\sigma x',x')\leq c_k^\frac{2}{n-2}C d_g(r_k\sigma x',2r_k\sigma x')\rightarrow Cd_{g_\infty}(\sigma x',2\sigma x')=\frac{C}{2\sigma}.
$$
The proof is completed  by applying triangle inequality.
\endproof

\begin{cor}\label{volume.density}
We have
$$
\lim_{\rho\rightarrow+\infty}
\frac{\vol(B_\rho^g(x_0),g)}{V_n\rho^n}=1
$$
and $(\overline{B}\setminus\{0\},\rho^{-2}g,x_0)$
converges to $(\R^n,0)$ in the Gromov-Hausdorff distance
for any $x_0$.
\end{cor}

\proof
First, we prove the convergence of volume ratio.
It suffices to prove that for any $\rho_k\rightarrow+\infty$,
a subsequence of $\frac{\vol(B_{\rho_k}^g(x_0),g)}{V_n{\rho_k}^n}$ converges to 1.

Let $\rho_k\rightarrow+\infty$, and $x_k=a_kx_0$ for some $a_k\in\R^+$, such that
$$
d_g(x_0,x_k)=\rho_k.
$$
Put $y_k=\sigma x_k$, where $\sigma=2^m$ is sufficiently large.
We denote
$$
\tau_k=d_g(x_k,y_k).
$$
We will first approximate $B_{\rho_k}^g(x_0)$ with $B_{\tau_k}^g(y_k)$, and subsequently approximate $B_{\tau_k}^g(y_k)$ by its intersection with $B_{|y_k|}$, that is, $B_{\tau_k}^g(y_k) \cap B_{|y_k|}$. The reason for doing this is that after rescaling, 
$B_{\tau_k}^g(y_k) \cap B_{|y_k|}$ exhibits very good convergence properties.\\

By \eqref{distance.decay},
$$
d_g(x_0,y_k)<3d_g(y_k,\frac{1}{2}y_k):=\sigma_k,
$$
hence
\begin{equation}\label{tau.rho}
\tau_k-\sigma_k\leq\rho_k\leq\tau_k+\sigma_k.
\end{equation}

Since
$$
d_g(x,x_0)\leq d_g(x,y_k)+d_g(x_0,y_k),
$$
for any $x\in B_{\tau_k-2\sigma_k}^g(y_k)$,
$$
d_g(x,x_0)\leq \tau_k-2\sigma_k+\sigma_k\leq\rho_k.
$$
Similarly, for any $x\in B_{\rho_k}^g(x_0)$
$$
d_g(x,y_k)\leq d_g(x,x_0)+d_g(x_0,y_k)\leq\rho_k+\sigma_k\leq \tau_k+2\sigma_k.
$$
Then
\begin{equation}\label{GH001}
B^g_{\tau_k-2\sigma_k}(y_k)\subset B^g_{\rho_k}(x_0),\s
B^g_{\rho_k}(x_0)\subset B^g_{\tau_k+2\sigma_k}(y_k).
\end{equation}
Let $u_k=u(r_kx)r^\frac{n-2}{2}c_k$, where $r_k=|x_k|$ and $c_k$ is as in Lemma \ref{asy}. By Lemma \ref{asy},
$u_k$ converges to $|x|^{2-n}$, and
\begin{equation}\label{sigma.tau}
\frac{\sigma_k}{\tau_k}\rightarrow \frac{3}{\sigma-1}.
\end{equation}
Then, by choosing $\sigma$ sufficiently large, for a fixed $\epsilon$ and
sufficiently large $k$,
\begin{equation}\label{BxBy}
B^g_{(1-\epsilon)\tau_k}(y_k)\subset B^g_{\rho_k}(x_0)\subset B^g_{(1+\epsilon)\tau_k}(y_k).
\end{equation}

Next, we estimate
$$
\frac{\vol(B_{\lambda\tau_k}^g(y_k),g)}{(\lambda\tau_k)^n}.
$$
By \eqref{volume.estimate}, \eqref{volume.over.distance} and \eqref{distance.decay}, we have
\begin{eqnarray*}
\frac{\vol(B_{\lambda\tau_k}^g(y_k)\cap B_{{|y_k|}},g)}{{(\lambda\tau_k)}^n}&\leq&
\frac{\vol(B_{{\lambda\tau_k}}^g(y_k),g)}{{(\lambda\tau_k)}^n}\\
&=&\frac{\vol(B_{\lambda \tau_k}^g(y_k)\cap B_{|y_k|},g)+
\vol(B_{\lambda \tau_k}^g(y_k)\setminus B_{|y_k|},g)}{(\lambda \tau_k)^n}\\
&\leq&\frac{\vol(B_{\lambda \tau_k}^g(y_k)\cap B_{{|y_k|}},g)+
\vol(B\setminus B_{{|y_k|}},g)}{(\lambda \tau_k)^n}\\
&\leq&\frac{\vol(B_{\lambda \tau_k}^g(y_k)\cap B_{{|y_k|}},g)+
C\vol(B_{2{|y_k|}}\setminus B_{{|y_k|}},g)}{(\lambda \tau_k)^n}\\
&\leq&\frac{\vol(B_{\lambda \tau_k}^g(y_k)\cap B_{{|y_k|}},g)}{(\lambda \tau_k)^n}+
C\frac{d_g^n(2y_k,y_k)}{(\lambda \tau_k)^n}\\
&\leq&\frac{\vol(B_{\lambda \tau_k}^g(y_k)\cap B_{{|y_k|}},g)}{(\lambda \tau_k)^n}+C
\frac{\sigma_k^n}{\lambda^n\tau_k^n}.
\end{eqnarray*}

It is easy to check that
$$
\frac{\vol(B_{\lambda\tau_k}^g(y_k)\cap B_{{|y_k|}},g)}{(\lambda\tau_k)^n}\rightarrow\frac{\vol(B^{g_\infty}_{\lambda(1-\sigma^{-1}) }(\sigma x_0/|x_0|)\cap B_{\sigma },g_\infty)}{(\lambda(1-\sigma^{-1}))^n}.
$$
To calculate $\vol(B^{g_\infty}_{(1-\sigma^{-1}) }(\sigma x_0/|x_0|)\cap B_{\sigma},g_\infty)$, we use the coordinates change
$$
B\rightarrow \R^n\setminus B: \s x\rightarrow z=\frac{x}{|x|^2}.
$$
In the new coordinates, $g_\infty=g_{euc}$, and
$$
B^{g_\infty}_{\lambda(1-\sigma^{-1}) }(\sigma x_0/|x_0|)\cap B_{\sigma }=
B_{\lambda(1-\sigma^{-1})}(\sigma^{-1}x_0/|x_0|)\setminus B_{\sigma^{-1}}.
$$
Thus
$$
\vol(B^{g_\infty}_{\lambda(1-\sigma^{-1}) }(\sigma x_0/|x_0|)\cap B_{\sigma},g_\infty)=V_n((\lambda(1-\sigma^{-1}))^n-O(\sigma^{-n})).
$$
Then we can select $\sigma$ to be sufficiently large and let $\lambda=1\pm\epsilon$,  such that
$$
1-\epsilon\leq\frac{\vol(B_{(1+\epsilon)\tau_k}^g(y_k),g)}{V_n((1+\epsilon)\tau_k)^n},\s \frac{\vol(B_{(1-\epsilon)\tau_k}^g(y_k),g)}{V_n((1-\epsilon)\tau_k)^n}\leq(1+\epsilon)
$$
when $k$ is sufficiently large. By \eqref{BxBy} and \eqref{tau.rho},
$$
1-C\epsilon\leq\frac{\vol(B_{\rho_k}^g(x_0))}{V_n\rho_k^n}\leq 1+C\epsilon
$$
when $k$ is sufficiently large, we complete the proof of the ratio convergence.\\

Next, we prove the Gromov-Hausdorff convergence. It suffices to prove that
a subsequence of $(\overline{B_1^{g/\rho_k^2}(x_0)},d_{g/\rho_k^2},x_0)$ converges to $(\overline{B},0)$.

By \eqref{tau.rho}, \eqref{GH001}, and \eqref{sigma.tau},
$$
\lim_{\sigma\rightarrow+\infty}\lim_{k\rightarrow+\infty}d_{GH}\left((\overline{B_1^{g/\rho_k^2}(x_0)},d_{g/\rho_k^2},x_0),(\overline{B_1^{g/\tau_k^2}(y_k)},d_{g/\tau_k^2},y_k)\right)=0.
$$
Combining \eqref{sphere.distance.decay} with \eqref{annulus.diam.decay}, we have
$$
diam(B_1\setminus B_{|y_k|},g)<Cd(\partial B_{|y_k|},\partial B_{|y_k|/2},g)\leq C\tau_k\frac{d(\partial B_{|y_k|},\partial B_{|y_k|/2},g)}{d(x_k,y_k)}\leq \frac{C}{\sigma}\tau_k,
$$
since $\frac{d(\partial B_{|y_k|},\partial B_{|y_k|/2},g)}{d(x_k,y_k)}
\rightarrow \frac{1/(2\sigma)}{1-1/\sigma}$.
Then
$$
\lim_{\sigma\rightarrow+\infty}\lim_{k\rightarrow+\infty}d_{GH}\left((\overline{B_1^{g/\tau_k^2}(y_k)},d_{g/\tau_k^2},y_k),(\overline{B_1^{g/\tau_k^2}(y_k)\cap B_{|y_k|}},d_{g/\tau_k^2},y_k)\right)=0.
$$
Note that
$$
({B_1^{g/\tau_k^2}(y_k)\cap B_{|y_k|}},d_{g/\tau_k^2},y_k)=({B_1^{g_k/(c_k^\frac{4}{n-2}\tau_k^2)}(y_k/r_k)\cap B_{|y_k/r_k|}},d_{g_k/(c_k^\frac{4}{n-2}\tau_k^2)},y_k/r_k).
$$
Since $d_{g_k}(x_k,y_k)\rightarrow 1-1/\sigma$ and $d_{g_k}=c_k^\frac{2}{n-2}d_g$, we have
$c_k^\frac{2}{n-2}\tau_k\rightarrow 1-1/\sigma$, which implies that
$$
\lim_{\sigma\rightarrow+\infty}\lim_{k\rightarrow+\infty}d_{GH}\left(\overline{B_1^{g/\tau_k^2}(y_k)\cap B_{|y_k|}},d_{g/\tau_k^2},y_k),(\overline{B_{1}^{g_k}(y_k/r_k)\cap B_{|y_k/r_k|}},d_{g_k},y_k/r_k)\right)=0.
$$
However, we have
$$
(\overline{B_{1}^{g_k}(y_k/r_k)\cap B_{|y_k/r_k|}},d_{g_k},y_k/r_k)\rightarrow (\overline{B_{1}^{g_\infty}(\sigma x_0/|x_0|)\cap B_{\sigma}},d_{g_\infty},\sigma x_0/|x_0|),
$$
and the limit is isometric to
$$
(\overline{B_{1}(\sigma^{-1}x_0/|x_0|)\setminus B_{\sigma^{-1}}},d_{g_{euc}},\sigma^{-1}x_0/|x_0|).
$$
Letting $\sigma\rightarrow \infty$, we complete the proof.
\endproof

\begin{center}
\begin{tikzpicture}



\draw[line width=0.7pt](1.55,0)--(9.5,0);
\filldraw[line width=0.5pt,opacity=1, pattern=dots,even odd rule]
 (9.5,0) circle (0.5) (9.5,0) circle (2.9);

\draw[line width=1.3pt] (9,0) circle (1.3pt);
\draw[line width=1.3pt] (8,0) circle (1.3pt);
\draw (9.55,0) circle (2pt);
\draw[line width=1.3pt] (5.7,0) circle (1.3pt);
\draw[dashed] (9.55,0) circle (1.55);
\filldraw[opacity=0.2,even odd rule]
 (9.5,0) circle (0.5) (9.55,0) circle (1.55);

\draw (9.8,0) node { $0$};
\draw (9.3,-0.3) node { $x_k$};
\draw (8,-0.3) node {$y_k$};
\draw (5.8,-0.3) node { $x_0$};
\draw (10,-1.1) node { $B_{|y_k|}$};

\draw (4,1.5) node { $\rho_k=d_{g}(x_0,x_k)$};
\draw (4,1) node { $\tau_k=d_{g}(x_k,y_k)$};
\draw (4.2,2) node { $y_k=\sigma x_k=2^mx_k$};

\clip (7.5,0) circle (1.5);
\end{tikzpicture}

{\tiny \it Fig 1.  $B_{\tau_k}^g(y_k)$ is the whole region filled with small dots. The shaded region is $B_{\tau_k}^g(y_k)\cap B_{|y_k|}$.\\
We use $B_{\tau_k}^{g}(y_k)\cap B_{|y_k|}$ to approximate $B_{\rho_k}^{g}(x_0)$.}
\end{center}\vspace{1ex}

\begin{cor}\label{rigidity}
Assume $(M,g)$ is conformally equivalent to a domain of a compact manifold without boundary.
If $\|R\|_{L^\frac{n}{2}}<+\infty$ and $Ric\geq 0$, and if $(M,g)$ is complete and noncompact, then $(M,g)=\R^n$.
\end{cor}

\proof By the Bishop-Gromov Theorem, $\vol(B_r^g(x),g)\leq V_nr^n$. Then, by a result in \cite{Carron-Herzlich}, there exists a compact manifold $(M_0,g_0)$ and a finite set $A\subset M_0$, such that $(M,g)$ is conformal to $(M_0\setminus A,g_0)$.
By Corollary \ref{volume.density}, $A$
contains a single point, so the corollary follows from the Bishop-Gromov Theorem.
\endproof

Next, we derive a stronger version of Lemma \ref{asy} and
finish the proof of Theorem \ref{main}:
\begin{pro}
$w=G^{-1}u$ is in $W^{2,p}(B)$ for any $p\in[1,\frac{n}{2})$, where $G$ is the Green function defined by
$$
-\Delta_{g_0}G=\delta_0,\s G|_{\partial B}=0.
$$
\end{pro}

\proof By direct computation,
$$
\Delta_{g_0} u=G\Delta_{g_0} w+2\nabla_{g_0} G\nabla_{g_0} w=-c(n)Ru^\frac{n+2}{n-2}+c(n)R(g_0)u.
$$
Then
\begin{eqnarray*}
-\Delta_{g_0} w&=&c(n)G^{-1}Ru^\frac{n+2}{n-2}+2\nabla_{g_0}\log G\nabla_{g_0} w-c(n)R(g_0)w\\
&=&c(n)Ru^\frac{4}{n-2}w+2\nabla_{g_0}\log G\nabla_{g_0} w-c(n)R(g_0)w\\
&:=&f.
\end{eqnarray*}
It is well known (cf.\cite{Aubin}) that $G=r^{2-n}(1+O(1))$ near 0 and
$$
|\nabla_{g_0}\log G|(x)\leq \frac{C}{|x|}
$$
when $x$ is small.

First, we show $w\in L^q$ for any $q$. Indeed, applying \eqref{u.decay} to $\alpha=q$, $\beta=(n-2)q$,
we get
$$
\int_{B_t}|w|^q\leq C\sum_i\int_{B_{2^{-i}t}}|x|^{q(n-2)}|u|^q<C\sum_i((1+\tau)2^{-n})^i<+\infty.
$$

Next, we show $w\in W^{1,p}(B)$ for any $p<n$. Since
$$
|\nabla_{g_0}w|\leq G^{-1}|\nabla_{g_0}u|+uG^{-1}|\nabla_{g_0}\log G|
\leq C(|x|^{n-2}|\nabla_{g_0}u|+|x|^{n-3}u),
$$
we may apply \eqref{u.decay} to $\alpha=p$ and $\beta=(n-3)p$ to
get
$$
\int_{B_t}(|x|^{n-3}|u|)^p<+\infty
$$
and \eqref{nabla.decay} to $\alpha=p$ and $\beta=(n-2)p$ to obtain
$$
\int_{B_t}|x|^{n-2}|\nabla u|^p<+\infty.
$$
Then $\int_{B_t}|\nabla w|^p<+\infty$ for any $p<n$. Let $\varphi\in\mathcal{D}(B)$, and $\eta_\epsilon$ be a cutoff function which is
1 in $B\setminus B_{2\epsilon}$, 0 in $B_\epsilon$, and satisfies $|\nabla\eta_\epsilon|<\frac{C}{\epsilon}$. It is easy to check that $\eta_\epsilon w$ is bounded in $W^{1,p}(B)$, hence a subsequence of $\eta_\epsilon w$
converges weakly in $W^{1,p}(B)$. Obviously, $w$ is the limit, hence $w\in
W^{1,p}(B)$.

Next, we show $f\in L^p$ for any $p<\frac{n}{2}$ and $w$ solves the equation $-\Delta_{g_0}w=f$ weakly in $B$.
Since
$$
|f|\leq C(|R|u^\frac{4}{n-2}w+|x|^{-1}|\nabla_{g_0}w|+w),
$$
by the fact that $|R|u^\frac{4}{n-2}\in L^\frac{n}{2}$, $w\in L^q$ for any $q>0$ and $\nabla w\in L^p$ for any $p<n$,
it is easy to check that $f\in L^p$ for any $p<\frac{n}{2}$. Then
$$
\int_{B}\nabla_{g_0}\varphi\nabla_{g_0} wdV_{g_0}=
\lim_{\epsilon\rightarrow 0}\int_B\nabla_{g_0}\eta_\epsilon\varphi\nabla_{g_0} wdV_{g_0}
=\lim_{\epsilon\rightarrow 0}\int_B\eta_\epsilon\varphi f=\int_B\varphi f,
$$
hence $w$ is a weak solution.

The proof can be completed without difficulty using the theory of elliptic equations.
\endproof

\section{conformally immersed submanifolds in $\R^{n+k}$}
In this section, we consider a conformal immersion $F: (B \setminus {0}, g_0) \rightarrow \R^{n+k}$ satisfying $|A|_{L^n} < +\infty$, where $A$ represents the second fundamental form.  We define
$$
g=F^*(g_{euc})=u^\frac{4}{n-2}g_0.
$$
Obviously
$$
\int_{B}|R|^\frac{n}{2}dV_g<+\infty.
$$
For the purposes of this section, we assume $\vol(F(B \setminus {0})) = +\infty$. As a consequence of Corollary \ref{complete.metric}, the space $(\overline{B}_{\frac{1}{2}} \setminus {0}, g)$ is complete. 
The goal of this section is to prove Theorem \ref{as.map}. \\

{\it Proof of Theorem \ref{as.map}:} First of all,
we can find $c_k$, such that $c_k r_k^\frac{n-2}{2} u(r_k x)$
converges to $|x|^{2-n}$ weakly in $W^{2,p}_{loc}(\R^n\setminus\{0\})$
for any $p<\frac{n}{2}$.
Set
$$
F_k=c_k^\frac{2}{n-2}(F(r_kx)-F(r_kx_0))+y_0,
$$
where $y_0$ will be defined later.
Since
$$
\left|\frac{\partial F_k}{\partial x^i}\right|^2=
r_k^2c_k^\frac{4}{n-2}(u(r_kx))^\frac{4}{n-2}=(c_kr_k^\frac{n-2}{2}u(r_kx))^\frac{4}{n-2},
$$
$|\nabla F_k|$ is bounded in $W^{2,p}(B_r\setminus B_\frac{1}{r})$ by Lemma \ref{asy}. Thus, we may assume
$F_k$ converges weakly in $W^{3,p}_{loc}(\R^n\setminus\{0\})$ to a map $F_\infty$ which satisfies
$F_\infty(\R^n\setminus\{0\})\subset \R^n$ and
$$
\frac{\partial F_\infty}{\partial x^i}\frac{\partial F_\infty}{\partial x^j}=|x|^{-4}\delta_{ij}.
$$
For convenience, we transition to new coordinates
$$
x\rightarrow y=\frac{x}{|x|^2}.
$$
In these coordinates,
$$
g_\infty(y)=g_{euc}(y).
$$
Then $F_\infty$ can be considered as an isometric map from $\R^n\setminus\{0\}$
to $\R^n$.

Let $\gamma(t)=ty$, where $y\in S^{n-1}$. Since $\gamma(t)$ is a geodesic
in $\R^n\setminus\{0\}$, $f(\gamma(t))$ must be a ray of $\R^n$.
In addition, it is easy to check that as $y'$, $y''$
approach 0
$d_{g_\infty}(F_\infty(y'),F_\infty(y''))\rightarrow 0$. Then $\lim\limits_{y\rightarrow 0}F_\infty(y)$ exists. Select $y_0$ such that the limit is $\lim\limits_{y\rightarrow 0}F_\infty(y)=0$. Then
$F_\infty(\gamma(t))=tF_\infty(y)$.

Since $F_\infty$ is isometric, for any $X\in T\partial B_1=S^{n-1}$ with $|X|=1$, it holds
$$
F_{\infty,*}(X)\bot F_{\infty,*}(\frac{\partial}{\partial r}),
\s |F_{\infty,*}(X)|=1.
$$
Hence, the restriction $F_\infty|_{ S^{n-1}}$ is an isometric map from
$S^{n-1}$ to itself. Since $S^{n-1}$ is simply connected, $F_\infty|_{S^{n-1}}$ is a homeomorphism, as follows from the fact that
an isometric map is a covering map. Therefore, we may assume $F_\infty(y)=y$.

In the original coordinates, this translates to
$$
F_\infty(x)=\frac{x}{|x|^2}.
$$
Proceeding, we consider a sequence $x_k \rightarrow 0$. Assuming $a = |x_0|$ and setting $x_k' = a\frac{x_k}{|x_k|}$, $r_k = |x_k|$ and $\sigma = 2^{m}$, we find that
$$
d_g(\sigma x_k,x_k')\leq Cd_g(\sigma x_k,2\sigma x_k),\s and \s \frac{
d_g(\sigma x_k,2\sigma x_k)}{d_g(x_k,\sigma x_k)}\rightarrow\frac{1}{2(\sigma-1)}.
$$
Thus we can choose $m$, such that for large $k$, we have
$$
1-\epsilon\leq \frac{d_g(x_k,x_0)}{d_g(x_k,\sigma x_k)}\leq 1+\epsilon.
$$
By  arguments as in the proof of Lemma \ref{decay}, for any $\tau\in(0,1)$,
$$
\frac{1-\tau}{2}\leq\frac{|F(rx)-F(rx/2)|}{|F(rx/2)-F(rx/4)|}\leq \frac{1+\tau}{2}
$$
for any $x\in\partial B$ and sufficiently small $r$. Then we can choose $m$, such that for sufficiently large $k$,
$$
1-\epsilon\leq\frac{|F(x_k)-F(x_0)|}{|F(x_k)-F(\sigma x_k)|}<1+\epsilon.
$$

Since
$$
\frac{|F(x_k)-F(\sigma x_k)|}{d_g(x_k,\sigma x_k)}=
\frac{|F_k(x_k/|x_k|)-F_k(\sigma x_k/|x_k|)|}{d_{g_k}(x_k/|x_k|,\sigma x_k/|x_k|)}\rightarrow 1,
$$
it follows that
$$
1-C\epsilon<\frac{|F(x_k)-F(x_0)|}{d_g(x_k,x_0)}\leq 1+C\epsilon
$$
when $k$ is sufficiently large.
\endproof

\section{ 4 dimensional Gauss-Bonnet-Chern formulas}
In this section, we assume $n=4$ and discuss Gauss-Bonnet-Chern formulas.

For our purpose, we set $dS_{g_0}=\Theta(r,\theta)dS^3$ and define $\phi=\log u$ and
$$
F_1(r)=\int_{S^3}R(r,\theta)u^{2}(r,\theta)\Theta(r,\theta)dS^3,\s H_1=
-\int_{S^3}R(r,\theta)\frac{\partial }{\partial r}(u^{2}(r,\theta)\Theta(r,\theta))dS^3
$$
$$
F_2(r)=\int_{S^3}\left(\Delta_{g_0}\phi(r,\theta)\right)\Theta(r,\theta)dS^3,\s H_2=-
\int_{S^3}\left(\Delta_{g_0}\phi\right)\frac{\partial }{\partial r}\Theta(r,\theta)dS^3.
$$
Let $n_{g,\partial B_r}$ be the unit normal vector of $\partial B_r$ with respect to $g$.
If we choose $x^1,\cdots,x^n$ to be normal coordinates of $g_0$, then $B_r=B_r^{g_0}(0)$, $
n_{g,\partial B_r}=u^{-1}\frac{\partial}{\partial r}$,and
\begin{equation}\label{boundary.term1}
\int_{\partial B_r}n_{g,\partial B_r}(R)dS_g=\int_{S^3}u^2\frac{\partial R}{\partial r}\Theta dS^3=
F_1'(r)+H_1(r).
\end{equation}
Moreover, we have
\begin{equation}\label{boundary.term2}
\int_{\partial B_r}\frac{\partial \Delta_{g_0}\phi}{\partial r}dS_{g_0}=
F_2'(r)+H_2(r).
\end{equation}
Equations \eqref{boundary.term1} and \eqref{boundary.term2} will help us to calculate $\int_{B_r}\Delta_gRdV_g$ and $\int_{B_r}\Delta_{g_0}^2\log udV_{g_0}$.
For example, by
$$
\int_{B_r}\Delta_gRdV_g=\lim_{r_k\rightarrow 0}(F_1'+H_1(r))|_{r_k}^{r},
$$
if we can find a sequence $r_k\rightarrow 0$, such that the limit of $F_i'(r_k)+H_i(r_k)$ is known, then we will get the exact value of $\int_{B_r}\Delta_gRdV_g$.

The following lemma will play a vital role in the following discussions.

\begin{lem}\label{mean.value}
Let $f\in C^1[\frac{r_0}{4}, 2r_0]$, $h\in C^0[\frac{r_0}{4}, 2r_0]$ and $b_1$, $b_2$ are constants. Assume
$$
\left|\int_{\frac{r_0}{4}}^{\frac{r_0}{2}}fdt-\frac{3}{32}b_1r_0^2\right|+\left|\int_{r_0}^{2r_0}fdt-\frac{3}{2}b_1r_0^2\right|\leq ar_0^2,\s \int_{\frac{r_0}{4}}^{2r_0}|h-b_2|<ar_0.
$$
Then there exists $\xi\in[r_0/4,2r_0]$ such that
$$
|f'(\xi)+h(\xi)-b_1-b_2|\leq 12a.
$$
\end{lem}

\proof Since we can replace $f$ with $f-b_1r$ and $h$ with $h-b_2$, it suffices for our aim to prove the case when $b_1=b_2=0$.
By The Mean Value Theorem for Integrals, there exists $\xi_1\in[r_0/4,r_0/2]$ and $\xi_2\in[r_0,2r_0]$, such that
$$
\frac{r_0}{4}f(\xi_1)=\int_\frac{r_0}{4}^\frac{r_0}{2}f(t)dt,\s
r_0f(\xi_2)=\int_{r_0}^{2r_0}f(t)dt,
$$
which yields that
$$
|f(\xi_1)|\leq 4ar_0,\s |f(\xi_2)|\leq ar_0.
$$
Then
$$
\left|\int_{\xi_1}^{\xi_2}(f'+h)\right|\leq |f(\xi_1)-f(\xi_2)|+\int_{\xi_1}^{\xi_2}|h|
\leq 5ar_0+\int_{r_0/4}^{2r_0}|h|\leq 6ar_0.
$$
Using the Mean Value Theorem for Integrals again, we can find $\xi\in[\xi_1,\xi_2]$, such that
$$
(\xi_2-\xi_1)|f'(\xi)+h(\xi)|\leq 6a r_0.
$$
Noting that $r_0/2<\xi_2-\xi_1$, we complete the proof.
\endproof

\begin{lem}
For any sufficiently small $r$, we have
$$
\left|\int_r^{2r}F_1(t)dt\right|<\alpha(r)r^2, \s \int_r^{2r}|H_1|(t)dt<\alpha(r)r,
$$
and
$$
\left|\int_r^{2r}F_2(t)dt+6\omega_3r^2\right|<\alpha(r)r^2, \s \int_r^{2r}|H_2-12\omega_3|(t)dt<\alpha(r)r,
$$
where $\lim\limits_{r\rightarrow 0}\alpha(r)=0$, and $\omega_3=2\pi^2$ is the volume of the $3$-dimensional sphere.
\end{lem}
\proof
We have
\begin{eqnarray*}
\left|\int_{r}^{2r} F_1(t)dt\right|
& \leq& \int_{r}^{2r} \int_{S^3} |R(g)| u^2 \Theta dS^3dt  = \int_{B_{2r}\setminus B_r}|R(g)| u^2dV_{g_0}\\
&\leq&
\left(
\int_{B_{2r} \setminus B_{r}} | R(g) u^2 |^2 dV_{g_0}
\right)^{\frac{1}{2}}
\left(
\int_{B_{2r} \setminus B_{r}} dV_{g_0}
\right)^{\frac{1}{2}} \\
&\leq& C\|R\|_{L^2(B_{2r},g)}r^2.
\end{eqnarray*}
and
\begin{eqnarray*}
\int_{r}^{2r} |H_1(t)| dt
& \leq& 2\int_{r}^{2r} \int_{S^3}
|R(g)| u^2 \left|\frac{\partial \log u}{\partial t}\right|\Theta dS^3dt
+ \int_{r}^{2r} \int_{S^3}
|R(g)| u^2\left|\frac{\partial\log \Theta(t,\theta)}{\partial t}\right|(\Theta dS^3dt) \\
&\leq& C(\|R\|_{L^2(B_{2r},g)}\|\nabla \phi\|_{L^2(B_{2r}\setminus B_r,g_0)}+\int_{B_{2r}\setminus B_r}|R(g)| u^2\frac{1}{r}dV_{g_0})\\
&\leq& C\|R\|_{L^2(B_{2r},g)}\|\nabla \phi\|_{L^2(B_{2r}\setminus B_r,g_0)}+Cr\|R(g)\|_{L^2(B_{2r},g)}.
\end{eqnarray*}
By Lemma \ref{Lalpha},
$$
\|\nabla \phi\|_{L^2(B_{2r}\setminus B_r,g_0)}\leq Cr\|Ru^2+R(g_0)\|_{L^2(B_{2r}\setminus B_r,g_0)}\leq C(\|R\|_{L^2(B_{2r}\setminus B_r,g)}r+r^2),
$$
hence
$$
\int_r^{2r}|H_1(t)|dt\leq C(\|R\|_{L^2(B_{2r}\setminus B_r,g)}r+r^2).
$$

Next, we discuss $F_2$. Since
$$
-\Delta_{g_0}\phi-|\nabla_{g_0}\phi|^2=c(n)Ru^2-c(n)R(g_0),
$$
we obtain
\begin{eqnarray*}
\int_{r}^{2r} F_2(t) dt
& =& \int_{B_{2r}\setminus B_r}\Delta_{g_0}\phi dV_{g_0} \\
& =& -\int_{B_{2r}\setminus B_r}|\nabla_{g_0} \phi|^2dV_{g_0}-c(4)\int_{B_{2r}\setminus B_r}(R(g) u^2-R(g_0))dV_{g_0}.
\end{eqnarray*}
Note that
$$
\int_{B_{2r}\setminus B_r}|R(g) u^2-R(g_0)|dV_{g_0}\leq C(\|R\|_{L^2(B_{2r}\setminus B_r,g)}r^2+r^4).
$$
To get the estimate of $\int_{r}^{2r}F_2$, we only need to show that
$$
\lim_{r\rightarrow 0}\frac{1}{r^2}\int_{B_{2r}\setminus B_r}|\nabla_{g_0}\phi|^2dV_{g_0}=6\omega_3.
$$

Assume there exists $r_k\rightarrow 0$, such that
$$
\lim_{k\rightarrow \infty}\frac{1}{r_k^2}\int_{B_{2r_k}\setminus B_{r_k}}|\nabla_{g_0}\phi|^2dV_{g_0}=\lambda\neq 6\omega_3.
$$
Set $u_k=c_kr_ku(r_kx)$, where $c_k$ is chosen such that $\int_{\partial B_1}\log u_k=0$.
By the arguments in Section 3, $\log u_k(x)$ converges to
$\log |x|^{-2}$ weakly in $W^{2,p}_{loc}(\R^4\setminus\{0\})$. Then, after passing to a subsequence,
$\int_{B_2\setminus B_1}|\nabla_{g_0(r_kx)}\log u_k|^2$ converges to $6\omega_3$, hence
$$
\lim_{k\rightarrow+\infty}\frac{1}{r_k^2}\int_{B_{2r_k}\setminus B_{r_k}}|\nabla_{g_0}\phi|^2dV_{g_0}=
\lim_{k\rightarrow+\infty}\int_{B_2\setminus B_1}|\nabla_{g_0(r_kx)}\log u_k|^2dV_{g_0(r_kx)}=6\omega_3.
$$
This leads to a contradiction.\\

Lastly, we calculate $\int_r^{2r}|H_2-12\omega_3|$:
\begin{equation*}
\int_{r}^{2r} |H_2(t) -12\omega_3|dt
\leq \int_{B_{2r}\setminus B_r}\left||\nabla_{g_0} \phi|^2\frac{\partial\log\Theta}{\partial t}-\frac{12}{\Theta}\right|dV_{g_0}+C\int_{B_{2r}\setminus B_r}|R(g) u^2-R(g_0)|\frac{1}{t}dV_{g_0}.
\end{equation*}
The same argument as above shows that
$$
\lim_{r\rightarrow 0}\frac{1}{r}\int_{B_{2r}\setminus B_r}\left||\nabla_{g_0} \phi|^2\frac{\partial\log\Theta}{\partial t}-\frac{12}{\Theta}\right|dV_{g_0}=0.
$$
This completes the proof.
\endproof

We will provide several applications here. First, we calculate $\int_{B_r}\Delta_{g} R$:
\begin{lem}
There exists $r_k \to 0$ such that
$$
\int_{B_{\frac{1}{2}}\setminus B_{r_k}}\Delta_g R(g) dV_g\rightarrow\int_{\partial B_\frac{1}{2}}
\frac{\partial R}{\partial r}dS_g.
$$
\end{lem}

\proof
We have
$$
\int_{B_{\frac{1}{2}}\setminus B_{r}}\Delta_g R(g) dV_g=\int_{\partial B_\frac{1}{2}}n_{g,\partial B_\frac{1}{2}}(R)dS_g-
\int_{\partial B_r}n_{g,\partial B_r}(R)dS_g.
$$
Applying Lemma \ref{mean.value} to $b_1=b_2=0$ and $f=F_1$, $h=H_1$ , we deduce this lemma from \eqref{boundary.term1}.
\endproof

We recall some basic properties of $Q$-curvatures, cf. \cite{Chang1,Chang2}. On a 4-dimensional manifold, the Paneitz operator is defined as follows:
$$
P_{g_0}\varphi=\Delta_{g_0}^2\varphi+div_{g_0}\left(\frac{2}{3}R_{g_0}\nabla_{g_0}\varphi-2Ric_{g_0}^{ij}\varphi_i\frac{\partial}{\partial x^j}\right).
$$
The $Q$-curvature of $g$ satisfies the following equations:
$$
Q(g)=-\frac{1}{12}\Delta_gR(g)-\frac{1}{4}|Ric(g)|^2+\frac{1}{12}R^2,
$$
$$
P_{g_0}\phi+2Q(g_0)=2Q(g)e^{4\phi}.
$$
For simplicity, we define
$$
T(\phi)=\frac{1}{3}R_{g_0}\frac{\partial \phi}{\partial r}-Ric_{g_0}(\nabla_{g_0}\phi,\frac{\partial}{\partial r}).
$$

\begin{lem} \label{GBC1.local}
There exists $r_k\rightarrow 0$, such that
$$
\lim_{k\rightarrow+\infty}\int_{B_\frac{1}{2}\setminus B_{r_k}} Q(g)dV_g=\int_{B_\frac{1}{2}}Q(g_0)dV_{g_0}-4\omega_3+
\int_{\partial B_\frac{1}{2}}\left(\frac{1}{2}\frac{\partial\Delta_{g_0}\phi}{\partial r}
+T(\phi)\right).
$$
\end{lem}

\proof
We have
\begin{eqnarray*}
\int_{B_\frac{1}{2}\setminus B_r}Q_gdV_g&=&\int_{B_\frac{1}{2}\setminus B_r}Q_{g_0}dV_{g_0}+\frac{1}{2}\int_{B_\frac{1}{2}\setminus B_r}
P_{g_0}(\phi)dV_{g_0}\\
&=&\int_{B_\frac{1}{2}\setminus B_r}Q_{g_0}dV_{g_0}+\frac{1}{2}\int_{\partial(B_\frac{1}{2}\setminus B_r)}\left(\frac{\partial\Delta_{g_0}\phi}{\partial r}
+\frac{2}{3}R_{g_0}\frac{\partial \phi}{\partial r}-2Ric_{g_0}(\nabla_{g_0}\phi,\frac{\partial}{\partial r})\right)dS_{g_0}.
\end{eqnarray*}
By \eqref{logu.decay},
$$
\lim_{r\rightarrow 0}\int_{\partial B_r}|\nabla_{g_0}\phi|dS_{g_0}= 0.
$$
Then
\begin{equation}\label{nabla.phi}
\int_{\partial B_r}\left(\frac{2}{3}R_{g_0}\frac{\partial \phi}{\partial r}-2Ric_{g_0}(\nabla_{g_0}\phi,\frac{\partial}{\partial r})\right)dS_{g_0}
\rightarrow 0.
\end{equation}

By applying Lemma \ref{mean.value} to $f=F_2/2$, $h=H_2/2$ and $(b_1,b_2)=(-2\omega_3,6\omega_3)$, we deduce from \eqref{boundary.term2}
that there exists $r_k$, such that
$$
\frac{1}{2}\int_{\partial B_{r_k}}\frac{\partial\Delta_{g_0}\phi}{\partial r}dS_{g_0}\rightarrow 4\omega_3.
$$
Therefore, we complete the proof.
\endproof

Next, we discuss the relationship between $\|R\|_{L^2}$
and $\|Riem\|_{L^2}$:
\begin{lem} \label{R2Ric2.local}
We have
\begin{equation*}
\int_{B_\frac{1}{2}}|Riem(g)|^2dV_g<+\infty.
\end{equation*}
\end{lem}

\proof
It is well-known that
$$
Riem(g)=W(g)+\frac{1}{2}(Ric(g)-\frac{1}{6}R(g)g){~\wedge\!\!\!\!\!\bigcirc~}g,
$$
where $W$ is the Weyl tensor and ${~\wedge\!\!\!\!\!\bigcirc~}$ is the Kulkarni-Nomizu product. Since $|W|^2dV_g$ is conformally invariant, we only need to check
$Ric(g)\in L^2$ here.

Recall that
$Q(g)=-\frac{1}{12}\Delta_gR(g)-\frac{1}{4}|Ric(g)|^2+\frac{1}{12}R^2$,
which means that
\begin{eqnarray*}
\int_{B_\frac{1}{2}\setminus B_r}|Ric(g)|^2dV_g&=&\frac{1}{3}\int_{B_\frac{1}{2}\setminus B_r}R^2dV_{g}-\frac{1}{3}\int_{B_\frac{1}{2}\setminus B_r}\Delta_gR(g)dV_{g}-4\int_{B_\frac{1}{2}\setminus B_r}Q(g)dV_{g}\\
&=&\frac{1}{3}\int_{B_\frac{1}{2}\setminus B_r}R^2dV_{g}-4\int_{B_\frac{1}{2}\setminus B_r}Q(g_0)dV_{g_0}\\
&&-\frac{1}{3}\int_{\partial (B_\frac{1}{2}\setminus B_r)}\frac{\partial R}{\partial r}dS_g-2\int_{\partial (B_\frac{1}{2}\setminus B_r)}(\frac{\partial \Delta_{g_0}\phi}{\partial r}+2T(\phi))dS_{g_0}
\end{eqnarray*}
Applying Lemma \ref{mean.value} to $f=\frac{1}{3}F_1+2F_2$, $h=\frac{1}{3}H_1+2H_2$, $(b_1,b_2)=(-8\omega_3,24\omega_3)$, we can find $r_k\rightarrow 0$, such that
$$
\frac{1}{3}\int_{\partial B_{r_k}}\frac{\partial R}{\partial r}dS_g+2\int_{\partial B_{r_k}}(\frac{\partial \Delta_{g_0}\phi}{\partial r}+2T(\phi))\rightarrow 16\omega_3.
$$
~\endproof

Lastly, we consider the formula for Pfaffian form:
\begin{lem}\label{GBC2.local}
We have
\begin{eqnarray*}
\int_{B_\frac{1}{2}}Pf(g)
&=&-4\omega_3+\int_{B_\frac{1}{2}}Pf(g_0)
+\frac{1}{2}\int_{\partial B_\frac{1}{2}} \frac{\partial \Delta_{g_0}\phi}{\partial r}dS_{g_0}\\
&&+\frac{1}{12}\int_{\partial B_\frac{1}{2}}u^2\frac{\partial R}{\partial r}dS_{g_0}-\frac{1}{12}\int_{B_\frac{1}{2}}\Delta_{g_0}R(g_0)dV_{g_0}+\int_{\partial B_\frac{1}{2}}T(\phi)dS_{g_0}.
\end{eqnarray*}
\end{lem}
\proof
Recall that
$$
Pf(g)=\frac{1}{8}|W(g)|^2+\frac{1}{12}R^2-\frac{1}{4}|Ric(g)|^2,
$$
where $W$ is the Weyl tensor. Since
$$
\int_{B}|W(g)|^2dV_{g}=\int_B|W(g_0)|^2dV_{g_0}<+\infty,
$$
$Pf(g)$ is integrable. Recall that (c.f. \cite{Chang2})
\begin{equation*}
\begin{aligned}
Pf(g) & = \frac{1}{8} |W(g)|^2 dV_g
+ Q(g) dV_g
+ \frac{1}{12} \Delta_g R(g) dV_g, \\
Pf(g_0) & = \frac{1}{8} |W(g_0)|^2 dV_{g_0}
+ Q(g_0) dV_{g_0}
+ \frac{1}{12} \Delta_{g_0} R(g_0) dV_{g_0},
\end{aligned}
\end{equation*}
and
\begin{equation*}
P_{g_0} \phi + 2Q(g_0) =2 Q(g) e^{4 \phi},
\end{equation*}
where $g = u^2 g_0 = e^{2 \phi} g_0$,
we have
\begin{eqnarray*}
\int_{B_\frac{1}{2}\setminus B_r} Pf(g)
& =& \int_{B_\frac{1}{2}\setminus B_r} Pf(g_0)
+ \frac{1}{2}\int_{B_\frac{1}{2}\setminus B_r} P_{g_0} \phi dV_{g_0}\\
&&	+ \frac{1}{12} \int_{B_\frac{1}{2}\setminus B_r} \Delta_g R(g) dV_g
-\frac{1}{12} \int_{B_\frac{1}{2}\setminus B_r} \Delta_{g_0} R(g_0) dV_{g_0}\\
& =&\int_{B_\frac{1}{2}\setminus B_r}Pf(g_0)
+ \frac{1}{2}\int_{\partial(B_\frac{1}{2}\setminus B_r)} \frac{\partial \Delta_{g_0}\phi}{\partial r}dS_{g_0}\\
&&+\frac{1}{12}\int_{\partial (B_\frac{1}{2}\setminus B_r)}u^2\frac{\partial R}{\partial r}dS_{g_0}-\frac{1}{12}\int_{B_\frac{1}{2}\setminus B_r}\Delta_{g_0}R(g_0)dV_{g_0}\\
&&+\frac{1}{2}\int_{\partial(B_\frac{1}{2}\setminus B_r)}\left(R_{g_0}\frac{\partial \phi}{\partial r}-2Ric_{g_0}(\nabla_{g_0}\phi,\frac{\partial}{\partial r})\right)dS_{g_0}.
\end{eqnarray*}
Apply Lemma \ref{mean.value} to $f=\frac{1}{12}F_1+\frac{1}{2}F_2$,
and $h=\frac{1}{12}H_1+\frac{1}{2}H_2$, $(b_1,b_2)=(-2\omega_3,6\omega_3)$, which suffices to complete the proof.
\endproof

Theorem \ref{GBC1}, \ref{GBC2} and \ref{R2Ric2} can be deduced from
Lemma \ref{GBC1.local}, \ref{GBC2.local} and \ref{R2Ric2.local} easily.
\\

\section{Examples}
In the last section,
we provide  examples of metrics on $B^4_{1/2}\setminus\{0\}$ that are conformal to $g_{euc}$ and satisfy  $\|R(g)\|_{L^2}<+\infty$. We will set
$$
u=r^{-2}e^{v},\s \phi=\log u,\s and \s g=u^2g_{euc},
$$
where $v=v(r)$ is radial.

We have
$$
|R_g|^2dV_g=(c\frac{\Delta u}{u})^2dx,\s |Q_g|dV_g=|\Delta^2\phi|dx,
$$
where $c$ is a constant. We observe that
\begin{equation}\label{exa.1}
\frac{\Delta u}{u}=v''-\frac{v'}{r}+(v')^2,\s 
|\Delta^2\phi|=|\Delta^2 v|\leq C(|v''''|+\frac{|v'''|}{r}+\frac{|v''|}{r^2}+\frac{|v'|}{r^3}).
\end{equation}

Recall that (cf. \cite[Ch. 5]{Schoen-Yau})
$$
R_{ij}(g)=-(\log u^2)_{,ij}+\frac{1}{2}(\log u^2)_i(\log u^2)_j-\frac{1}{2}(\Delta(\log u^2)+|\nabla\log u^2|^2)(g_{euc})_{ij}.
$$
Note that  the Hessian tensor of $\log u^2$ in the euclidean metric is
$$
Hess(\log u^2,g_{euc})=(\log u^2)''dr\otimes dr+r(\log u^2)'g_{S^3}.
$$
It follows that
\begin{eqnarray*}
Ric(g)&=&\left(-(\log u^2)''+\frac{1}{2}|(\log u^2)'|^2-\frac{1}{2}(\Delta(\log u^2)+|(\log u^2)'|^2)\right)dr\otimes dr\\
&&-\left(r(\log u^2)'+\frac{1}{2}(\Delta(\log u^2)+|(\log u^2)'|^2)r^2\right)g_{S^3}\\
&=&-\frac{3}{2}\left((\log u^2)''+\frac{1}{r}(\log u^2)'\right)dr\otimes dr\\
&&-\left(\frac{5}{2}r(\log u^2)'+\frac{r^2}{2}(\log u^2)''+\frac{r^2}{2}|(\log u^2)'|^2\right)g_{S^3}\\
&=&-3\left(v''+\frac{1}{r}v'\right)dr\otimes dr-\left(-3rv'+r^2v''+2|v'|^2r^2\right)g_{S^3},
\end{eqnarray*}
leading to
\begin{equation}\label{exa.2}
|Ric(g)|^2\sqrt{|g|}\leq C(|v''|^2+\frac{|v'|^2}{r^2}).
\end{equation}

\begin{eg}
Consider $v=r^a\log r$,  $g=e^{-2(2-r^a)\log r}g_{euc}$. 
We have
$$
\frac{\Delta u}{u}= r^{2a - 2}(a\log r+1)^2 + r^{a - 2}(a^2\log r - 2a\log r + 2a-2). 
$$
Then, $R(g)\in L^2$ if and only $a>0$. In this setting, it is easy to  verified that $Ric(g)\in L^2$
and $Q(g)\in L^1$ from \eqref{exa.1} and \eqref{exa.2}.\\
\end{eg}

\begin{eg}\label{Eg2}
Let $v=-a\log(-\log r)$, $g=\frac{g_{euc}}{r^4|\log r|^{2a}}$. We find
$$
\frac{\Delta u}{u}=\frac{a(1+a)}{r^2\log^2 r}+\frac{2a}{r^2\log r}.
$$
Then, $R(g)\in L^2$ for any $a$. We can check that $Ric(g)\in L^2$ and $Q(g)$ is integrable. 

This example extends the metric  
$$
g=\frac{|dz|^2}{|z|^2|\log|z||^{2a}}.
$$
constructed by Hulin-Troyanov \cite{Hulin-Troyanov} on a 2 dimensional disk, which has finite total Gauss curvature. Depending on the value of $a$, their metric can be either bounded or unbounded, finite area or infinite area. However, in our case, the metric is always unbounded and of infinite volume.

Note that $r^2u=e^{v}=|\log r|^{-a}$ does not belong to $W^{2,2}$ when $a>-\frac{1}{2}$. This indicates that the conclusion `$G^{-1}u\in W^{2,p}$ for any $p<\frac{n}{2}$' in Theorem \ref{main} can not be extended to $p=\frac{n}{2}$.
\end{eg}
 
\begin{eg}\label{Eg3}
Consider $v=r^4\sin\frac{1}{r}$, $g=r^{-4}e^{2r^4\sin\frac{1}{r}}g_{euc}$. We observe that
$$
v''=O(1),\s v'/r^2=O(1),\s \Delta^2 v=\frac{\sin\frac{1}{r}}{r^4}+O(\frac{1}{r^3}).
$$
Consequently, the scalar curvature $R(g)$ and the Ricci curvature $Ric(g)$ are in $L^2$. Since
$$
\int_{B_\frac{1}{2}}\frac{|\sin(\frac{1}{r})|}{r^4}dx=\int_0^\frac{1}{2}\frac{|\sin(\frac{1}{r})|}{r^4}r^3dr=\int_2^\infty\frac{|\sin(t)|}{t}dt>\sum_{k=2}^\infty\frac{1}{k\pi}\int_{(k-1)\pi}^{k\pi}|\sin(t)|dt=+\infty,
$$
the Q-curvature $Q(g)$ is not integral in this case.
\end{eg}

\end{document}